\documentclass[a4paper,11pt]{amsart}

\usepackage{fancybox}
\usepackage{amscd}
\usepackage{amsmath}
\usepackage{amssymb}
\usepackage{amsthm}
\usepackage{float}
\usepackage[dvips, pdftex]{graphicx}
\usepackage{tikz}
\usepackage[all]{xy}

\usepackage{color}

\usepackage[mathscr]{eucal}
\usepackage{stmaryrd}
\usepackage[hidelinks, hyperfootnotes=false]{hyperref}
\usepackage{lipsum}

\newcommand\blfootnote[1]{%
  \begingroup
  \renewcommand\thefootnote{}\footnote{#1}%
  \addtocounter{footnote}{-1}%
  \endgroup
}

\author{Michail Savvas}
\title[Stabilizer reduction and DTK invariants]{Intrinsic Stabilizer Reduction and Generalized Donaldson--Thomas Invariants}
\address{Department of Mathematics, The University of Texas at Austin, Austin, TX 78712, USA}
\email{msavvas@utexas.edu}

\usepackage{array}
\usepackage{amscd}

\DeclareFontFamily{U}{rsfs}{%
\skewchar\font127}
\DeclareFontShape{U}{rsfs}{m}{n}{%
<-6>rsfs5<6-8.5>rsfs7<8.5->rsfs10}{}
\DeclareSymbolFont{rsfs}{U}{rsfs}{m}{n}
\DeclareSymbolFontAlphabet
{\mathrsfs}{rsfs}
\DeclareRobustCommand*\rsfs{%
\@fontswitch\relax\mathrsfs}

\theoremstyle{plain}
\newtheorem{theorem}{Theorem}[section]
\newtheorem*{theorem*}{Theorem}
\newtheorem{proposition}[theorem]{Proposition}
\newtheorem{lemma}[theorem]{Lemma}

\newtheorem{definition}[theorem]{Definition}
\newtheorem{remark}[theorem]{Remark}
\newtheorem{corollary}[theorem]{Corollary}

\newtheorem{prop-defi}[theorem]{Proposition-Definition}
\newtheorem{thm-defi}[theorem]{Theorem-Definition}
\newtheorem{thm-constr}[theorem]{Theorem-Construction}
\newtheorem{lem-defi}[theorem]{Lemma-Definition}

\newtheorem{exam}[theorem]{Example}

\newtheorem{setup}[theorem]{Setup}

\newtheorem{setup-def}[theorem]{Setup-Definition}

\newdimen\argwidth
\def\db[#1\db]{
 \setbox0=\hbox{$#1$}\argwidth=\wd0
 \setbox0=\hbox{$\left[\box0\right]$}
  \advance\argwidth by -\wd0
 \left[\kern.3\argwidth\box0 \kern.3\argwidth\right]}

\newcommand{\aA}{\mathcal{A}}

\newcommand{\cC}{\mathcal{C}}

\newcommand{\fF}{\mathcal{F}}

\newcommand{\kK}{\mathcal{K}}

\newcommand{\mM}{\mathcal{M}}
\newcommand{\nN}{\mathcal{N}}
\newcommand{\oO}{\mathcal{O}}
\newcommand{\pP}{\mathcal{P}}

\newcommand{\sS}{\mathcal{S}}

\newcommand{\uU}{\mathcal{U}}

\newcommand{\yY}{\mathcal{Y}}
\newcommand{\zZ}{\mathcal{Z}}

\newcommand{\Ob}{\mathcal{O}b}
\newcommand{\OB}{\mathop{\rm Ob}\nolimits}
\newcommand{\Hone}{\mathop{H^1}\nolimits}
\newcommand{\bL}{\mathbb{L}}
\newcommand{\bQ}{\mathbb{Q}}

\newcommand{\fm}{\mathfrak{m}}
\newcommand{\fg}{\mathfrak{g}}
\newcommand{\fr}{\mathfrak{r}}
\newcommand{\fh}{\mathfrak{h}}
\newcommand{\fc}{\mathfrak{c}}

\renewcommand{\tilde}{\widetilde}
\renewcommand{\hat}{\widehat}

\newcommand{\tU}{\tilde{U}}

\newcommand{\lr}{\longrightarrow}

\newcommand{\tX}{\widetilde{X}}

\newcommand{\id}{\textrm{id}}

\newcommand{\rk}{\mathop{\rm rk}\nolimits}

\newcommand{\Ext}{\mathop{\rm Ext}\nolimits}
\newcommand{\Spec}{\mathop{\rm Spec}\nolimits}

\newcommand{\Coh}{\mathop{\rm Coh}\nolimits}
\newcommand{\QCoh}{\mathop{\rm QCoh}\nolimits}

\newcommand{\Perf}{\mathop{\rm Perf}\nolimits}

\newcommand{\im}{\mathop{\rm im}\nolimits}

\newcommand{\Aut}{\mathop{\rm Aut}\nolimits}

\newcommand{\Crit}{\mathop{\rm Crit}\nolimits}

\newcommand{\GL}{\mathop{\rm GL}\nolimits}

\newcommand{\bA}{\mathbb{A}}
\newcommand{\bC}{\mathbb{C}}
\newcommand{\bZ}{\mathbb{Z}}
\newcommand{\bP}{\mathbb{P}}
\newcommand{\bG}{\mathbb{G}}
\newcommand{\bT}{\mathbb{T}}

\def\lal{_\lambda}
\def\vir{\mathrm{\vir}}

\def\lra{\longrightarrow}

\def\lalp{_\alpha}
\def\sub{\subset}

\def\virt{^{\mathrm{vir}}}

\def\uss{^{ss}}
\def\beq{\begin{equation}}
\def\eeq{\end{equation}}

\def\git{/\!\!/}
\def\lalp{_\alpha }
\def\lab{_{\alpha\beta}}
\def\lbet{_\beta}

\def\hV{\hat V}
\def\hU{\hat U}
\def\uss{^{ss}}

\newcommand{\cdga}{{\mathrm{cdga_{\bC}^{\leq 0}}}}
\newcommand{\dSpec}{\mathbf{Spec}\,}
\newcommand{\dSt}{\mathbf{dSt}_{\bC}}
\newcommand{\qcoh}{\mathsf{qcoh}}
\newcommand{\dgmod}{\mathrm{dg}\textrm{-}\mathrm{Mod}}
\newcommand{\dR}{\mathrm{dR}}
\newcommand{\boldit}[1]{\boldsymbol{#1}}
\DeclareMathAlphabet\BCal{OMS}{cmsy}{b}{n}
\newcommand{\dd}{\mathrm{d}}
\newcommand{\derived}{\mathbf{D}}

\makeatletter
 
  \@addtoreset{equation}{section}
\makeatother

\makeatletter
\def\@tocline#1#2#3#4#5#6#7{\relax
  \ifnum #1>\c@tocdepth 
  \else
    \par \addpenalty\@secpenalty\addvspace{#2}%
    \begingroup \hyphenpenalty\@M
    \@ifempty{#4}{%
      \@tempdima\csname r@tocindent\number#1\endcsname\relax
    }{%
      \@tempdima#4\relax
    }%
    \parindent\z@ \leftskip#3\relax \advance\leftskip\@tempdima\relax
    \rightskip\@pnumwidth plus4em \parfillskip-\@pnumwidth
    #5\leavevmode\hskip-\@tempdima
      \ifcase #1
       \or\or \hskip 1em \or \hskip 2em \else \hskip 3em \fi%
      #6\nobreak\relax
    \hfill\hbox to\@pnumwidth{\@tocpagenum{#7}}\par
    \nobreak
    \endgroup
  \fi}
\makeatother

\keywords{Donaldson--Thomas invariants; perfect complexes; Calabi--Yau
threefolds; intrinsic stabilizer reduction}
\subjclass{14A30, 14D20, 14D23, 14J30, 14J32, 14N35}

\begin{document}

\maketitle

\begin{abstract}
Let $\sigma$ be a stability condition on the bounded derived category $D^b({\mathop{\rm Coh}\nolimits} W)$ of a Calabi--Yau threefold $W$ and $\mathcal{M}$ a moduli stack parametrizing $\sigma$-semistable objects of fixed topological type. We define generalized Donaldson--Thomas invariants which act as virtual counts of objects in $\mathcal{M}$, fully generalizing the approach introduced by Kiem, Li and the author in the case of semistable sheaves. 

We construct an associated proper Deligne--Mumford stack $\widetilde{\mathcal{M}}^{\bC^\ast}$, called the $\mathbb{C}^\ast$-rigidified intrinsic stabilizer reduction of $\mathcal{M}$, with an induced semi-perfect obstruction theory of virtual dimension zero, and define the generalized Donaldson--Thomas invariant via Kirwan blowups to be the degree of the associated virtual cycle $[\widetilde{\mathcal{M}}]^{\mathrm{vir}} \in A_0(\widetilde{\mathcal{M}})$. This stays invariant under deformations of the complex structure of $W$. Applications include Bridgeland stability, polynomial stability, Gieseker and slope stability.
\end{abstract}

\tableofcontents

\section{Introduction}

\subsection{Brief historical background} Donaldson--Thomas (abbreviated as DT from now on) invariants constitute one of the main approaches for curve counting on Calabi--Yau threefolds. They naturally appear in many enumerative problems of interest in algebraic geometry and string theory and are conjecturally equivalent with other counting invariants, such as Gromov-Witten invariants, Stable Pair invariants \cite{PT1} and Gopakumar-Vafa invariants \cite{MaulikToda}. These relations have now been proven in many cases (for example, in \cite{MNOP1, MNOP2, TodaDTPT}).\blfootnote{The author was partially supported by a Stanford Graduate Fellowship, an Alexander S. Onassis Foundation Graduate Scholarship and an A.G. Leventis Foundation Grant during the course of this work.}

Let $W$ be a smooth, projective Calabi--Yau threefold and $\gamma \in H^*(W, \bQ)$. Classical DT theory was introduced in \cite{Thomas} in order to obtain virtual counts of stable sheaves on $W$ of Chern character $\gamma$. More precisely, Thomas considered the (coarse) moduli space $M := M^{ss}_L(\gamma)$ parameterizing Gieseker semistable sheaves on $W$ of positive rank, fixed determinant $L$ and Chern character $\gamma$. Assuming that every semistable sheaf is stable, $M$ is proper and admits a perfect obstruction theory in the sense of Li-Tian \cite{LiTian} or Behrend--Fantechi \cite{BehFan} of virtual dimension zero. Thus, there exists a virtual fundamental cycle $[M]\virt \in A_0 \left( M \right)$ and the classical DT invariant is defined as
$$\mathrm{DT}(M) := \deg \ [M]\virt.$$
One of its important properties is its invariance under deformation of the complex structure of $W$.
\medskip

Another important feature of DT invariants is their motivic nature. Behrend \cite{BehFun} showed that the obstruction theory of $M$ is symmetric in a certain sense and established an equality
\begin{align*} \label{Behrend identity} \tag{$\dagger$}
\mathrm{DT}(M) = \deg \ [M]\virt = \chi (M, \nu_{M} ) = \sum_{n \in \bZ} n \cdot \chi_{\mathrm{top}}(\nu_{M}^{-1}(n)),
\end{align*}
where $\nu_{M} \colon M \to \bZ$ is a canonical constructible function on $M$ and the right-hand side a weighted Euler characteristic.
\medskip

However, when stability and semistability of sheaves do not coincide, the above methods do not suffice and new techniques are required in order to define generalized DT invariants counting semistable sheaves. A main obstacle is that strictly semistable sheaves can have more automorphisms beyond $\bC^\ast$-scaling and, as a result, the stacks $\mM^{ss}(\gamma)$ parameterizing semistable sheaves are generally Artin and no longer Deligne--Mumford (after rigidifying $\bC^*$-scaling, cf. Subsection~\ref{Subsection 5.3}). This is necessary for the standard machinery of perfect obstruction theory and virtual cycles to apply.

In \cite{JoyceSong}, Joyce and Song constructed generalized DT invariants of moduli stacks $\mM^{ss}(\gamma)$ taking advantage of the above motivic behaviour and using motivic Hall algebras to obtain a generalization of the right-hand side of \eqref{Behrend identity}. Their DT invariant is easy to work with and amenable to computation.  
However, the proof of its deformation invariance is indirect and proceeds via wall-crossing to a stable pairs theory where semistability and stability coincide and thus a virtual cycle exists. Kontsevich and Soibelman \cite{KontSoibel} have also defined motivic generalized DT invariants using similar ideas. We also mention the related work of Behrend and Ronagh \cite{BehRon1, BehRon2}.
\medskip

In \cite{KLS}, the authors  develop a new direct approach towards defining generalized DT invariants of such a moduli stack $\mM^{ss}(\gamma)$. Their method adapts Kirwan's partial desingularization procedure \cite{Kirwan} to define an invariant as the degree of a zero-dimensional virtual cycle in an associated Deligne--Mumford stack $\tilde{\mM}$. The constructed invariant is called the generalized DT invariant via Kirwan blowups (or DTK invariant for short) and is a direct generalization of the left-hand side of \eqref{Behrend identity}. By the usual properties of virtual cycles, they establish the deformation invariance of DTK invariants.
\medskip

We make a final comment on the need for orientation data in the above approaches. Since motivic DT invariants are more refined, Kontsevich and Soibelman do need to assume the existence of such data to define their invariants. This assumption has been recently proved to hold in \cite{JoyUpOrient}. Orientation data are not necessary for the definition of numerical generalized DT invariants given in \cite{JoyceSong} and \cite{KLS}. This will be the case for the DTK invariant defined in this paper as well.

\subsection{Statement of results} This paper serves as a sequel to \cite{KLS}. We generalize the construction of DTK invariants to the case of moduli stacks $\mM := \mM^{\sigma-ss}(\gamma)$ parametrizing $\sigma$-semistable objects of Chern character $\gamma$ in the derived category $D^b(\Coh W)$ of coherent sheaves on $W$, where $\sigma$ is one of the following stability conditions:
\begin{enumerate}
\item A Bridgeland stability condition \cite{BridgStab}, as considered by Piyaratne--Toda \cite{TodaPiya} and \cite{quinticstab}.
\item A polynomial stability condition \cite{Bayer}, as considered by Lo \cite{Lo, Lo2}.
\item Gieseker and slope stability of sheaves. These are examples of weak stability conditions on $\Coh W$ in the sense of Joyce--Song \cite{JoyceSong}.
\end{enumerate}

More generally, $\sigma$ can be any nice stability condition (cf. Definition~\ref{stab cond} and Definition~\ref{nice stab cond}).

Our main results are summarized in the following theorem, giving a definition of DTK invariants counting $\sigma$-semistable complexes.

\begin{theorem*}
Let $W$ be a smooth, projective Calabi--Yau threefold, $\sigma$ a stability condition on $D^b(\Coh W)$ as in Definition~\ref{nice stab cond} and $\mM := \mM^{\sigma-ss}(\gamma)$ be the $(\bC^\ast$-rigidified$)$ moduli stack parametrizing $\sigma$-semistable complexes of Chern character $\gamma$. Then there exist:
\smallskip
\begin{enumerate}
    \item \emph{[Theorem-Construction~\ref{intrinsic stab red of stack with gms}, Theorem-Definition~\ref{rigidified intr stab red}]}\\ A canonical proper Deligne--Mumford stack $\tilde{\mM}^{\bC^\ast} \to \mM$, called the $\bC^\ast$-rigidified intrinsic stabilizer reduction of $\mM$. $\tilde{\mM}^{\bC^\ast}$ is isomorphic to $\mM$ over the open $\sigma$-stable locus ${\mM}^{\sigma-s}(\gamma) \subseteq \mM$.\smallskip
    \item \emph{[Theorem~\ref{dt type admits vir cyc}, Theorem-Definition~\ref{DTK invariant of complexes}]}\\ A natural semi-perfect obstruction theory of virtual dimension zero on $\tilde{\mM}^{\bC^\ast}$, which extends the symmetric obstruction theory of the Deligne--Mumford stack ${\mM}^{\sigma-s}(\gamma)$. 
\end{enumerate}
\smallskip
We thus have a virtual fundamental cycle 
$$[\tilde{\mM}^{\bC^\ast}]\virt \in A_0(\tilde{\mM}^{\bC^\ast})$$ 
and the generalized Donaldson--Thomas invariant via Kirwan blowups of $\mM$ is defined as
$$\mathrm{DTK}(\mM^{\sigma-ss}(\gamma)) := \deg \ [\tilde{\mM}^{\bC^\ast}]\virt \in \bQ.$$
By Theorem~\ref{def inv of dtk}, this is invariant under deforming the complex structure of $W$.
\end{theorem*}

\subsection{Brief review of the case of sheaves and sketch of construction} We first give a brief account of the results of \cite{KLS} and then explain the necessary adjustments in order to generalize their approach.
\medskip

Let $\mM$ be an Artin stack parametrizing Gieseker semistable sheaves of a fixed Chern character on a smooth, projective Calabi--Yau threefold. By construction (see \cite{HuyLehn}), $\mM$ is obtained by Geometric Invariant Theory (GIT) (see \cite{MFK}), meaning that it is a quotient stack of the form
\begin{align} \label{GIT presentation}
    \mM = [X / G],
\end{align}
where $G$ is a reductive group acting on a projective space $\bP^N$ via a homomorphism $G\to GL(N+1,\bC)$ and $X$ is an invariant closed subscheme of the GIT semistable locus $X \subseteq (\bP^N)^{ss}$. 
\medskip

In order to define an invariant out of $\mM$, the first step in \cite{KLS} is the construction of a proper DM quotient stack $\tilde{\mM} := [\tX / G]$ over $\mM$, called the intrinsic stabilizer reduction of $\mM$. This is produced by an iterative blowup procedure, which successively resolves the loci of closed $G$-orbits in $X$ with stabilizer groups of maximum dimension. This procedure is an adaptation of Kirwan partial desingularization construction \cite{Kirwan} to singular GIT quotient stacks, using the notion of intrinsic blowup introduced in \cite{KiemLi} and suitably generalized in \cite{KLS}.
\medskip

The second step is to endow $\tilde{\mM}$ with a semi-perfect obstruction theory \cite{LiChang} of virtual dimension zero. This is a usual perfect obstruction theory in the sense of Behrend--Fantechi \cite{BehFan} defined on an \'{e}tale cover of $\tilde{\mM}$ together with compatibility data which play the role of suitable descent data. It is at this step where the existence of a virtual structure on $\mM$ plays a crucial role. More precisely, by \cite{PTVV}, $\mM$ is the truncation of a $(-1)$-shifted symplectic derived Artin stack and it then follows that $\mM$ is a d-critical stack \cite{JoyceDCrit}. In particular, applying Luna's \'{e}tale slice theorem \cite{Drezet} and using the theory of d-critical loci developed in \cite{JoyceDCrit}, we can obtain an explicit local description of $\mM$: for every closed point $x \in \mM$ with (reductive) stabilizer $H$, there exists a smooth affine scheme $V$ with an $H$-action, an invariant function $f \colon V \to \bA^1$ and an \'{e}tale morphism
\begin{align} \label{eq1}
[ U / H ] \to \mM,
\end{align}
containing $x$ in its image, where $U = \lbrace df = 0 \rbrace \sub V$ is the scheme-theoretic critical locus of $f$. Any two such local presentations can be compared as they give the same d-critical structure on $\mM$ in the sense of \cite{JoyceDCrit}.
\medskip 

To obtain the semi-perfect obstruction theory of $\tilde{\mM}$, consider the following $H$-equivariant 4-term complex on $U$:
\begin{equation}\label{eq2}
\fh=\mathrm{Lie}(H)\lra T_{V}|_{U} \xrightarrow{d(df)} F_{V}|_{U}=\Omega_{V}|_{U} \lra \fh^\vee.
\end{equation}
Around $u \in U$ with finite stabilizer group, this complex is quasi-isomorphic to a 2-term complex which gives a symmetric perfect obstruction theory of $[U / H]$ and thus of $\mM$ near $u$.

The intrinsic stabilizer reduction algorithm produces lifts of the \'{e}tale morphisms \eqref{eq1} to give an \'{e}tale cover
\begin{align} \label{eq3}
[ T / H ] \to \tilde{\mM},
\end{align}
where $T = \lbrace \omega_S = 0 \rbrace \subseteq S$ for $S$ a smooth affine $H$-scheme and $\omega_S \in H^0(S, F_S)$ an invariant section of an $H$-equivariant vector bundle $F_S$ on $S$. Moreover, there exists an effective invariant divisor $D_S$ such that \eqref{eq2} lifts to a 4-term complex
\begin{equation}\label{eq4}
\fh = \mathrm{Lie}(H) \lra T_{{S}}|_{{T}} \xrightarrow{d\omega_S} F_{{S}}|_{{T}} \lra \fh^\vee(-D_{{S}})
\end{equation}
whose first arrow is injective with locally free cokernel and last arrow is surjective. Therefore, \eqref{eq4} is quasi-isomorphic to a 2-term complex, whose dual can be shown to be a perfect obstruction theory on $[T/H]$. 
\medskip

In \cite{KLS}, it is shown that these perfect obstruction theories satisfy the axioms of a semi-perfect obstruction theory on $\tilde{\mM}$ of virtual dimension zero.

Since, as usual, semi-perfect obstruction theories produce virtual cycles, $\tilde{\mM}$ admits a virtual cycle $[\tilde{\mM}]\virt \in A_0(\tilde{\mM})$ whose degree defines the DTK invariant counting Gieseker semistable sheaves parametrized by $\mM$. A relative version of the above construction, using derived symplectic geometry, implies its deformation invariance.\\

In the present paper, we work with moduli stacks $\mM$ of semistable perfect complexes. These are truncations of $(-1)$-shifted symplectic derived Artin stacks as before, however they are not global quotient stacks obtained by GIT of the form \eqref{GIT presentation}. 
\medskip

Our main contribution is to remove this assumption. To do this, we use several recent major technical results on the \'{e}tale local structure of stacks \cite{AHR2, Alper}, moduli spaces of objects in abelian categories \cite{AlpHalpHein} and stability conditions in families \cite{familystab}.

We generalize the requirement of a GIT presentation to the condition that $\mM$ admit a good moduli space morphism \cite{GoodAlper}. We show that such stacks (under additional reasonable assumptions) admit an intrinsic stabilizer reduction $\tilde{\mM}$, which is a Deligne--Mumford stack and of independent interest in its own right (see Subsection~\ref{relations to derived ag} below).

We then proceed to show that if in addition $\mM$ is the truncation of a $(-1)$-shifted symplectic derived Artin stack, then the above statements about local models of $\mM$ and their comparison data carry over, using the main structural result of \cite{Alper} in place of Luna's \'{e}tale slice theorem. This allows us to define a semi-perfect obstruction theory on $\tilde{\mM}$ and the associated DTK invariant.
\medskip

Our approach thus works for a wide class of Artin stacks, which we call stacks of DT type (see Definition~\ref{stack of dt type}). Based on our discussion here, their two main characteristics can be summarized as the existence of a good moduli space and a derived enhancement that is $(-1)$-shifted symplectic. 

Moduli stacks of semistable complexes fit into this context. This is the case by the recent results of \cite{AlpHalpHein}, which combine the theory of good moduli spaces with the notion of $\Theta$-reductivity introduced in \cite{ThetaRed}, and apply to the stacks considered in the present paper.

Finally, the deformation invariance of the DTK invariant follows from the compatibility of the intrinsic stabilizer reduction with base change and the recent results on stability conditions in \cite{familystab}.

\subsection{Relation to derived algebraic geometry} \label{relations to derived ag} At the time of writing of this paper, intrinsic blowups and the intrinsic stabilizer reduction procedure were expected to be the classical shadow of a corresponding construction in derived algebraic geometry. 

Subsequently, derived blowups were first defined in full generality in the work of Hekking \cite{Hekking} and this expectation has recently been materialized by the results of \cite{HRS}, also using some of our results. Namely, the intrinsic blowup notion used in the present paper coincides with the classical truncation of an appropriate derived blowup of a derived scheme with a $G$-action along its derived fixed locus by the group $G$. Therefore, the intrinsic stabilizer reduction construction is indeed the classical truncation of a derived version of the Kirwan partial desingularization algorithm in GIT \cite{Kirwan} and its generalization to stacks with good moduli spaces by Edidin--Rydh \cite{EdidinRydh}, called derived stabilizer reduction. Moreover, the semi-perfect obstruction theory constructed here is closely related to the derived cotangent complex of the stabilizer reduction of a $(-1)$-shifted symplectic derived stack. 

\subsection{Layout of the paper} In \S \ref{sec background d-crit and sss} we review background material on d-critical loci and $(-1)$-shifted symplectic structures. In \S \ref{sec stacks}, we define the notion of stacks of DT type and establish their main properties that will be used throughout. \S \ref{sec desingularization} reviews the intrinsic stabilizer reduction procedure for GIT quotient stacks and generalizes it to stacks with good moduli spaces. In \S \ref{sec obstruction theory}, we explain how to construct a semi-perfect obstruction theory on stacks of DT type by following the arguments of \cite{KLS}. Finally, in \S \ref{sec DTK invariants}, we construct DTK invariants of semistable complexes by combining the above with the results of \cite{AlpHalpHein} and $\bG_m$-rigidification for Artin stacks and discuss their deformation invariance using recent work of \cite{familystab}.

\subsection{Acknowledgements} The author would like to thank his advisor Jun Li for introducing him to the subject, his constant encouragement and many enlightening discussions during the course of the completion of this work. He also benefitted greatly by conversations with Jarod Alper, Jack Hall, Daniel Halpern-Leistner, Young-Hoon Kiem, Alex Perry, David Rydh and Ravi Vakil.

\subsection{Notation and conventions} Here are the various notations and other conventions that we use throughout the paper:
\begin{itemize}
\item[--] All schemes and stacks are defined over the field of complex numbers $\bC$ or a smooth $\bC$-scheme $C$, unless stated otherwise. For this reason, reductive group schemes will be linearly reductive automatically. 

\item[--] $\mM$ typically denotes an Artin stack, of finite type, with affine stabilizers and separated diagonal unless stated otherwise.

\item[--] $W$ denotes a smooth, projective Calabi--Yau threefold over $\bC$ and $D^b(\Coh W)$ its bounded derived category of coherent sheaves.

\item[--] $C$ denotes a smooth quasi-projective scheme over $\bC$.

\item[--] $G$, $H$ denote complex reductive groups. Usually, $H$ will be a subgroup of $G$. $T$ denotes the torus $\bC^\ast$.

\item[--] If $x \in \mM$, $G_x$ denotes the automorphism group/stabilizer of $x$. We will only consider stabilizers of closed points $x \in \mM$. These will be reductive for most stacks of interest in this paper.

\item[--] If $U \hookrightarrow V$ is a closed embedding, $I_{U \subseteq V}$ or $I_U$ (when $V$ is clear from context) denotes the ideal sheaf of $U$ in $V$.

\item[--] For a morphism $\rho \colon U \to V$ and a sheaf $E$ on $V$, we systematically use $E \vert_U$ to denote $\rho^* E$, suppressing the pullback from the notation.

\item[--] If $V$ is a $G$-scheme, $V^G$ is used to denote the fixed point locus of $G$ in $V$.

\item[--] If $U$ is a scheme with a $G$-action, then $\hU$ is used to denote the Kirwan blowup of $U$ with respect to $G$. $\tU$ denotes the intrinsic stabilizer reduction of $U$.

\item[--] The abbreviations DT, DM, GIT, whenever used, stand for Donaldson--Thomas, Deligne--Mumford and Geometric Invariant Theory respectively.
\end{itemize}

\section{D-Critical Loci and $(-1)$-Shifted Symplectic Derived Stacks} \label{sec background d-crit and sss}

This section collects background material and terminology that will be used throughout the rest of the paper.

We first briefly recall Joyce's theory of d-critical loci, as developed in \cite{JoyceDCrit}, and establish some notation, and then proceed to quickly review shifted symplectic structures on derived stacks.

\subsection{d-critical schemes} \label{d-crit section} We begin by defining the notion of a d-critical chart.

\begin{definition}\emph{(d-critical chart)} \label{d-crit chart}
A d-critical chart for a scheme $X$ is the data of $(U,V,f,i)$ such that: $U \subseteq X$ is Zariski open, $V$ is a smooth scheme, $f \colon V \to \bA^1$ is a regular function on $V$ and $U \xrightarrow{i} V$ is a closed embedding so that $U = \lbrace d f=0 \rbrace = \Crit(f) \subseteq V$ is the scheme-theoretic vanishing locus of the derivative of $f$.

If $x \in U$, then we say that the d-critical chart $(U,V,f,i)$ is centered at $x$.
\end{definition}

Joyce defines a canonical sheaf $\sS_X$ of $\bC$-vector spaces with the property that for any Zariski open $U \subseteq X$ and an embedding $U \hookrightarrow V$ into a smooth scheme $V$ with ideal $I$, $\sS_X$ fits into an exact sequence
\begin{align} \label{loc 2.1}
0 \lr \sS_X |_U \lr \oO_V / I^2 \stackrel{d}{\lr} \Omega_{V} / I \cdot \Omega_{V}
\end{align}
For example, for a d-critical chart $(U,V,f,i)$ the element $f+I^2 \in \Gamma(V, \oO_V/I^2)$ gives a section of $\sS_X|_U$.

\begin{definition}\emph{(d-critical scheme)} \label{d-crit scheme}
A d-critical structure on a scheme $X$ is a section $s \in \Gamma(X, \sS_X)$ such that $X$ admits a cover by d-critical charts $(U,V,f,i)$ and $s|_U$ is given by $f+I^2$ as above on each such chart. We refer to the pair $(X,s)$ as a d-critical scheme.
\end{definition}

\subsection{Equivariant d-critical loci} For our purposes, we need equivariant analogues of the results of Subsection~\ref{d-crit section}. The theory works in parallel as before (cf. \cite[Section~2.6]{JoyceDCrit}).

\begin{definition} \emph{(Good action)} Let $G$ be an algebraic group acting on a scheme $X$. We say that the action is good if $X$ has a cover $\lbrace U_\alpha \rbrace_{\alpha \in A}$ where every $U_\alpha \subseteq X$ is a $G$-invariant Zariski open affine subscheme.
\end{definition}

\begin{remark}
If $X$ is affine, then trivially every action of $G$ on $M$ is good. 
\end{remark}

Suppose now that a scheme $X$ admits a good action by a reductive group $G$. Definitions~\ref{d-crit chart}, \ref{d-crit scheme} extend naturally to this equivariant setting (cf. \cite[Definition~2.40]{JoyceDCrit}) as follows: A $G$-invariant d-critical chart is given by the data $(U,V,f,i)$ of Definition~\ref{d-crit chart}, where $U \subseteq X$ is a $G$-invariant Zariski open subscheme, $V$ is a smooth scheme with a $G$-action, $i \colon U \to V$ an equivariant closed embedding and $f \colon V \to \bC$ a $G$-invariant function. The sheaf $\sS_X$ inherits a $G$-action by using a cover of $X$ by $G$-invariant open subschemes $U \subseteq X$ in \eqref{loc 2.1}. Then we say that the d-critical scheme $(X,s)$ is $G$-invariant if $s \in \Gamma(X, \sS_X)^G$ is a $G$-invariant global section of the sheaf $\sS_X$.

\begin{proposition} \emph{\cite[Remark~2.47]{JoyceDCrit}} \label{equivariant d-crit prop}
Let $G$ be a complex reductive group with a good action on a scheme $X$. Suppose that $(X,s)$ is an invariant d-critical scheme. Then for any $x \in X$ fixed by $G$, there exists an invariant d-critical chart $(U,V,f,i)$ centered at $x$, i.e. an invariant open affine $U \ni x$, a smooth scheme $V$ with a $G$-action, an invariant regular function $f \colon V \to \bA^1$ and an equivariant embedding $i \colon U \to V$ so that $U = \Crit(f) \subseteq V$.
\end{proposition}

\begin{remark} \label{Rmk d-crit}
If $G$ is a torus $\left( \bC^{\ast} \right)^k$, then Proposition \ref{equivariant d-crit prop} is true without the assumption that $x$ is a fixed point of $G$.
\end{remark}

\begin{remark}
One may replace Zariski open morphisms by \'{e}tale morphisms without any difference to the essence of the theory. Another option is to work in the complex analytic topology.
\end{remark}

\subsection{d-critical Artin stacks} The theory of d-critical loci extends naturally to Artin stacks. For more details, we point the interested reader to Section 2.8 of \cite{JoyceDCrit}. We mention the following definition and basic properties which we will need in the form of remarks.

\begin{definition} \emph{\cite[Corollary~2.52, Definition~2.53]{JoyceDCrit}} \label{d-crit artin}
Let $\mM$ be an Artin stack. For every smooth morphism $\phi \colon U \to \mM$, where $U$ is a scheme, the assignment $\sS(U, \phi) := \sS_U$ defines a sheaf $\sS_\mM$ of $\bC$-vector spaces in the lisse-\'{e}tale topology of $\mM$.

A d-critical structure on $\mM$ is a global section $s \in H^0(\mM, \sS_\mM)$. We then say that $(\mM, s)$ is a d-critical Artin stack.
\end{definition}

\begin{remark} \emph{\cite[Example~2.55]{JoyceDCrit}} \label{quot d-crit}
A d-critical structure can be equivalently described in terms of a groupoid presentation of $\mM$. It follows that d-critical structures on quotient stacks $\left[ X / G \right]$ are in bijective correspondence with invariant d-critical structures on $X$.

If $\mM' \to \mM$ is a smooth morphism of Artin stacks and $\mM$ is d-critical, then one may pull back the d-critical structure on $\mM$ making $\mM'$ d-critical as well.
\end{remark}

\subsection{Shifted symplectic structures}

Let $\cdga$ be the category of non-positively  graded commutative differential graded $\bC$-algebras. 

There is a spectrum functor $\dSpec \colon \cdga \to \dSt$ to the category of derived stacks (see \cite[Definition 2.2.2.14]{HAG-DAG} or \cite[Definition 4.2]{Toen_higher-derived}). An object of the form $\dSpec A$ is called an \emph{affine derived $\bC$-scheme}. Such objects provide the Zariski local charts for general \emph{derived} $\bC$-\emph{schemes}, see \cite[Section 4.2]{Toen_higher-derived}. An object $\BCal{M}$ in $\dSt$ is called a \emph{derived Artin stack} if it is $m$-geometric (cf.~\cite[Definition 1.3.3.1]{Toen_higher-derived}) for some $m$ and its `classical truncation' $t_0(\BCal{M})$ is an Artin stack (and not a higher stack). A derived Artin stack $\BCal{M}$ admits an \emph{atlas}, i.e.~a smooth surjective morphism $\boldit{U} \to \BCal{M}$ from a derived scheme. For a derived Artin stack $\BCal{M}$, there exists a cotangent complex $\bL_{\BCal{M}}$ of finite cohomological amplitude in $[-m,1]$ and a dual tangent complex $\bT_{\BCal{M}}$. Both are objects in a suitable stable $\infty$-category $L_{\qcoh}(\BCal{M})$ (see \cite{Toen_higher-derived} or \cite{HAG-DAG} for its definition).

Shifted symplectic structures on derived Artin stacks were introduced by Pantev--To\"{e}n--Vaqui\'e--Vezzosi in \cite{PTVV}. The definition is given in the affine case first, and then generalized by showing the local notion satisfies smooth descent. We recall the local definition: let us set $\boldit{M} = \dSpec A$, so that $L_{\qcoh}(\boldit{M}) \cong \derived(\dgmod_A)$. For all $p\geq 0$ one can define the exterior power complex $(\Lambda^p \bL_{\boldit{M}},\dd) \in L_{\qcoh}(\boldit{M})$, where the differential $\dd$ is induced by the differential of the algebra $A$. For a fixed $k \in \bZ$, define a $k$-shifted $p$-form on $\boldit{M}$ to be an element $\omega^0 \in (\Lambda^p \bL_{\boldit{M}})^k$ such that $\dd \omega^0 = 0$. To define the notion of closedness, consider the de Rham differential $\dd_{\dR}\colon \Lambda^p \bL_{\boldit{M}} \to \Lambda^{p+1} \bL_{\boldit{M}}$. A $k$-\emph{shifted closed $p$-form} is a sequence $(\omega^0,\omega^1,\ldots)$, with $\omega^i \in (\Lambda^{p+i}\bL_{\boldit{M}})^{k-i}$, such that $\dd \omega^0 = 0$ and $\dd_{\dR} \omega^i + \dd \omega^{i+1} = 0$. When $p=2$, any $k$-shifted $2$-form $\omega^0 \in (\Lambda^2 \bL_{\boldit{M}})^k$ induces a morphism $\omega^0\colon \bT_{\boldit{M}} \to \bL_{\boldit{M}}[k]$ in $L_{\qcoh}(\boldit{M})$, and we say that $\omega^0$ is \emph{non-degenerate} if this morphism is an isomorphism in $L_{\qcoh}(\boldit{M})$.

\begin{definition}[{\cite[Definition 1.18]{PTVV}}]
A $k$-shifted closed $2$-form $\omega = (\omega^0,\omega^1,\ldots)$ is called a $k$-\emph{shifted symplectic structure} if $\omega^0$ is non-degenerate. We say that $(\boldit{M},\omega)$ is a $k$-shifted symplectic (affine) derived scheme.
\end{definition}

When $k=-1$, it is shown in \cite{JoyceArt} that if $\mM$ is the classical truncation of a $(-1)$-shifted symplectic derived Artin stack, then $\mM$ admits an induced d-critical structure.

\begin{theorem}[{\cite[Theorem 3.18]{JoyceArt}}]\label{sympltoDcrit}
Let $(\BCal{M},\omega)$ be a $(-1)$-shifted symplectic derived Artin stack. Then the underlying classical Artin stack $\mM = t_0(\BCal{M})$ extends in a canonical way to a d-critical Artin stack $(\mM,s)$. This defines a `truncation functor' $\tau$ from the $\infty$-category of $(-1)$-shifted symplectic derived Artin stacks to the $2$-category of d-critical Artin stacks.
\end{theorem}

\section{Stacks of DT Type} \label{sec stacks}

In this section, we first give an account of results regarding the \'{e}tale local structure of Artin stacks and the theory of good moduli spaces. We then proceed to define stacks of DT type, develop standard local models for them and describe how to compare these models.

\subsection{Local structure of Artin stacks} The following theorem is an \'{e}tale slice theorem for stacks, which generalizes Luna's \'{e}tale slice theorem \cite{Drezet}. It states that Artin stacks are \'{e}tale locally quotient stacks.
\begin{theorem} \emph{\cite[Theorem~1.2]{Alper}} \label{Alper1}
Let $\mM$ be a quasi-separated Artin stack, locally of finite type over $\bC$ with affine stabilizers. Let $x \in \mM$ and $H \subseteq G_x$ a maximal reductive subgroup of the stabilizer of $x$. Then there exists an affine scheme $U$ with an action of $H$, a point $u \in U$ fixed by $H$, and a smooth morphism 
\begin{align*}
\Phi_x : \left[ U / H \right] \rightarrow \mM
\end{align*}
which maps $u \mapsto x$ and induces the inclusion $H \hookrightarrow G_x$ of stabilizers at $u$. Moreover, if $G_x$ is reductive, the morphism is \'{e}tale. If $\mM$ has affine diagonal, then $\Phi_x$ can be taken to be affine.
\end{theorem}

\subsection{Good moduli spaces} We now collect some useful results about the structure of a certain class of Artin stacks, namely those with affine diagonal admitting a good moduli space, following the theory developed by Alper \'{e}t al. All the material of the section can be found in \cite{GoodAlper} and \cite{Alper}.

We have the following definition of a good moduli space.
\begin{definition} \emph{\cite[Definition~4.1]{GoodAlper}} 
A morphism $q \colon \mM \to Y$, where $\mM$ is an Artin stack and $Y$ an algebraic space, is a good moduli space for $\mM$ if the following hold:
\begin{enumerate}
\item $q$ is quasi-compact and $q_* \colon \QCoh(\mM) \to \QCoh(Y)$ is exact.
\item The natural map $\oO_Y \to q_* \oO_\mM$ is an isomorphism.
\end{enumerate}
\end{definition}
The intuition behind the introduction of the notion of good moduli space is that stacks $\mM$ that admit good moduli spaces behave like quotient stacks $[X^{ss}/G]$ obtained from Geometric Invariant Theory (GIT) with good moduli space morphism given by the map $[X^{ss}/G] \to X \git G$. In this sense, it is a generalization of GIT quotients for stacks.

We state the following properties of stacks with good moduli space.
\begin{proposition} \emph{\cite[Proposition~4.7, Lemma~4.14, Theorem~4.16, Proposition 9.1, Proposition~12.14]{GoodAlper}} \label{properties of good mod spaces}
Let $\mM$ be locally noetherian over $\bC$ and $q \colon \mM \to Y$ be a good moduli space. Then:
\begin{enumerate}
\item $q$ is surjective.
\item $q$ is universally closed.
\item Two geometric points $x_1, x_2 \in \mM(k)$ are identified in $Y$ if and only if their closures $\overline{\lbrace x_1 \rbrace}$ and $\overline{\lbrace x_2 \rbrace}$ in $\mM \times_{\Spec \bC} \Spec k$ intersect.
\item Every closed point of $\mM$ has reductive stabilizer.
\item Let $y \in |Y|$ be a closed point. Then there exists a unique closed point $x \in |q^{-1}(y)|$.
\item If $\zZ$ is a closed substack of $\mM$, then the morphism $\zZ \to q(\zZ) = \im (\zZ)$ is a good moduli space morphism.
\item Let
$$\xymatrix{
\mM' \ar[r] \ar[d] & \mM \ar[d] \\
Y' \ar[r] & Y
}$$
be a Cartesian diagram of Artin stacks, such that $Y, Y'$ are algebraic spaces. 
\begin{enumerate}
\item If $\mM \to Y$ is a good moduli space, then $\mM' \to Y'$ is a good moduli space.
\item If $Y' \to Y$ is fpqc and $\mM' \to Y'$ is a good moduli space, then $\mM \to Y$ is a good moduli space.
\end{enumerate}

\item If $\mM$ is of finite type, then $Y$ is of finite type.
\end{enumerate}
\end{proposition}

The following notions will be useful for us.

\begin{definition} \emph{\cite[Definition~6.1, Remark~6.2]{GoodAlper}} \label{definition of saturated open substack}
Let $\mM$ be a locally noetherian Artin stack with a good moduli space $q \colon \mM \to M$. A Zariski open substack $\uU \subseteq \mM$ is called saturated if $q^{-1}(q(\uU)) = \uU$. 

If $\uU$ is saturated, the morphism $\uU \to U := q(\uU)$ gives a good moduli space for $\uU$ fitting into a Cartesian diagram
\begin{align*}
    \xymatrix{
    \uU \ar[d] \ar[r] & \mM \ar[d]^-{q} \\
    U \ar[r] & M,}
\end{align*}
where $U \to M$ is an open embedding.
\end{definition}

\begin{definition} \emph{\cite[Definition~3.13]{AHR2}} \label{definition of strongly etale}
Let $\mM',\ \mM$ be Artin stacks with good moduli space morphisms $q' \colon \mM' \to M'$, $q \colon \mM \to M$. A morphism $\mM' \to \mM$ is called strongly \'{e}tale if the induced morphism $M' \to M$ is \'{e}tale and the diagram
\begin{align*}
    \xymatrix{
    \mM' \ar[r] \ar[d]_-{q'} & \mM \ar[d]^-{q} \\
    M' \ar[r] & M
    }
\end{align*}
is Cartesian.
\end{definition}

\begin{remark}
In the sequel, we will be using the fact that strongly \'{e}tale morphisms are stabilizer-preserving.
\end{remark}

We have the following theorem regarding the \'{e}tale local structure of good moduli space morphisms for stacks with affine diagonal.

\begin{thm-defi} \emph{(Quotient chart) \cite[Theorem~4.12]{Alper}} \label{quotient chart existence thm-def}
Let $\mM$ be a locally noetherian Artin stack with a good moduli space $q \colon \mM \to M$ such that $q$ is of finite type with affine diagonal. If $x \in \mM$ is a closed point with (reductive) stabilizer $G_x$, then there exists an affine scheme $U$ with an action of $G_x$ and a Cartesian diagram
\begin{align} \label{strong quotient chart}
\xymatrix{
\left[ U_x / G_x \right] \ar[r]^-{\Phi_x} \ar[d] & \mM \ar[d]^-q \\
U_x \git G_x \ar[r] & M
}
\end{align}
such that $\Phi_x$ has the same properties as in Theorem~\ref{Alper1}, is affine and $U_x \git G_x$ is an \'{e}tale neighbourhood of $q(x)$.

We refer to any choice of data $(U_x, \Phi_x)$ such that $\Phi_x \colon [U_x / G_x] \to \mM$ is \'{e}tale, affine and stabilizer-preserving as a quotient chart for $\mM$ centered at $x$. We say that the quotient chart is strongly \'{e}tale if the morphism $\Phi_x$ is strongly \'{e}tale.
\end{thm-defi}

\subsection{Stacks of DT type} We start with the following definition.

\begin{definition} \emph{(Stack of DT type)} \label{stack of dt type}
Let $\mM$ be an Artin stack. We say that $\mM$ is of DT type if the following are true:
\begin{enumerate}
\item $\mM$ is quasi-separated and finite type over $\bC$.
\item There exists a good moduli space $q \colon \mM \to M$, where $q$ is of finite type and has affine diagonal.
\item $\mM$ is the classical truncation of a $(-1)$-shifted symplectic derived Artin stack $\BCal{M}$.
\end{enumerate}
\end{definition}

\begin{remark}
By Theorem~\ref{sympltoDcrit} a stack $\mM$ of DT type admits a canonical d-critical structure $s \in H^0(\mM, \sS_\mM)$.
\end{remark}

\begin{definition} \emph{(d-critical quotient chart)}
Let $x \in \mM$ and $\Phi \colon [U_x / G_x] \to \mM$ be as in Theorem-Definition~\ref{quotient chart existence thm-def}. We say that $(U_x, V, f, \Phi_x)$ is a d-critical quotient chart centered at $x$ if there exists a  $G_x$-invariant d-critical locus determined by data $(U_x, V, f, i)$.
\end{definition}

A variant of the following proposition first appeared and was used in \cite{TodaHall}.

\begin{proposition}\label{existence of d-crit quot}
Let $\mM$ be a stack of DT type and $x \in \mM$ a closed point. Then there exists a d-critical quotient chart for $\mM$ centered at $x$, which can be taken to be strongly \'{e}tale.
\end{proposition} 
\begin{proof} By Theorem-Definition~\ref{quotient chart existence thm-def}, we have a quotient chart $\Phi_x \colon [ U_x'' / G_x ] \to \mM$. By Proposition~\ref{equivariant d-crit prop}, Definition~\ref{d-crit artin} and Remark~\ref{quot d-crit}, there exists a $G_x$-invariant open affine $x \in U_x' \subseteq U_x''$ and an invariant d-critical chart $(U_x', V', f, i)$. We obtain a commutative diagram
\begin{align*}
    \xymatrix{
    [U_x' / G_x] \ar[r]^-{\Phi_x} \ar[d]_-{q'} & \mM \ar[d] \\
    U_x' \git G_x \ar[r] & M.
    }
\end{align*}
Applying the Fundamental Lemma \cite[Theorem~6.10]{AlperFundLemma} at $x$, there exists a Zariski open subscheme $R \subseteq U_x' \git G_x$ such that the induced diagram
\begin{align*}
     \xymatrix{
    (q')^{-1} (R) \ar[r]^-{\Phi_x} \ar[d]_-{q'} & \mM \ar[d] \\
    R \ar[r] & M
    }
\end{align*}
is Cartesian and $\Phi_x$ is strongly \'{e}tale. 

Since $U'_x \git G_x$ is affine, we may take $R$ to be affine as well. It is straightforward to check then that $(q')^{-1}(R) = [U_x / G_x]$, where $U_x$ is a $G_x$-invariant Zariski open affine subscheme of $U'_x$, containing $x$, and also $R = U_x \git G_x$.

Setting $V = V' \setminus (U' \setminus U_x)$, we have the d-critical chart $(U_x,V,f|_V,i)$. This concludes the proof.
\end{proof}

We now obtain the following key lemma, which gives a way to compare two choices of d-critical quotient charts.

\begin{lemma} \label{comparison of d-crit quot}
Let $\mM$ be a stack of DT type. Let $$(U\lalp, V\lalp, f\lalp, \Phi\lalp), \ (U\lbet, V\lbet, f\lbet, \Phi\lbet)$$ be two d-critical quotient charts with $V\lalp, V\lbet$ affine and $z \in [U\lalp / G\lalp] \times_\mM [U\lbet/G\lbet]$ a closed point with stabilizer $G := G_z$. Let $x\lalp \in [U\lalp/G\lalp]$ and $x\lbet \in [U\lbet/G\lbet]$ be the two projections of $z$. Then the following hold:
\begin{enumerate}
\item We have $G$-equivariant commutative diagrams for $\lambda = \alpha, \beta$
\begin{align} \label{loc 3.2}
\xymatrix{
T\lab \ar[r] \ar[d]_-{i\lal} & S\lab \ar[d]^-{\theta\lal} \\
U\lal \ar[r] & V\lal,
}
\end{align}
where the vertical arrows are unramified morphisms, the horizontal arrows are embeddings, $t \in T\lab$ maps to $x\lal \in U\lal$,  and $T\lab, S\lab$ are affine.
\item We have an induced diagram with \'{e}tale arrows:
\begin{align} \label{loc 3.3}
\xymatrix{
 & [T\lab/G] \ar[dr] \ar[dl] & \\
[U\lalp / G\lalp ] \ar[dr] & & [ U\lbet / G\lbet ] \ar[dl] \\
 & \mM &
 }
\end{align}
\item There exists a d-critical $G$-invariant quotient chart $(T\lab,S\lab,f\lab,\Phi\lab)$ for $\mM$ with $t \in T\lab$ fixed by $G$.
\item $I := (\theta\lal^* df\lal) = (df\lab)$ as ideal sheaves in $\oO_{S\lab}$ and $\theta\lal^* f\lal + I^2$, $f\lab + I^2$ give the same invariant section $s|_{T\lab} \in \Gamma (T\lab, \sS_{T\lab})^G$, where $\theta\lal$ are the morphisms defined above in diagram~\eqref{loc 3.2}.
\end{enumerate}
\end{lemma}

\begin{proof}
It follows by the second condition in Definition~\ref{stack of dt type} that $\mM$ must have affine diagonal. In particular, the Cartesian diagram
\begin{align*}
    \xymatrix{
    U\lalp \times_{\mM} U\lbet \ar[d] \ar[r] & \mM \ar[d]^-{\Delta_{\mM}} \\
    U\lalp \times U\lbet \ar[r] & \mM \times \mM}
\end{align*}
implies that $U\lab := U\lalp \times_{\mM} U\lbet$ is affine.

Let us write $G\lab = G\lalp \times G\lbet$. Since $U\lab$ and $V\lalp \times V\lbet$ are affine and the product embedding $U\lalp \times U\lbet \hookrightarrow V\lalp \times V\lbet$ is $G\lab$-invariant, there exists a smooth affine scheme $V\lab$ and a $G\lab$-equivariant embedding fitting into a commutative diagram
\begin{align*}
    \xymatrix{
    U\lab \ar[d] \ar[r] & V\lab \ar[d] \\
    U\lalp \times U\lbet \ar[r] & V\lalp \times V\lbet.}
\end{align*}

Concretely, we may choose any $G\lab$-equivariant embedding $U\lab \hookrightarrow \bA^N$ for big enough $N$ and then take $V\lab = V\lalp \times V\lbet \times \bA^N$ with the obvious choices of morphisms that make the above diagram commute.

Let $t \in U\lab$ a representative of the point $z \in [U\lalp / G\lalp] \times_\mM [U\lbet/G\lbet] \simeq [U\lab / G\lab]$ (we are liberally abusing notation at this point). Using Luna's \'{e}tale slice theorem, take $S\lab \subseteq V\lab$ to be an affine \'{e}tale slice for $t$ in $V\lab$ and $T\lab := S\lab \cap U\lab$ the induced affine \'{e}tale slice for $t \in U\lab$.

It is clear that conditions (1) and (2) are now satisfied and remain so even after $G$-invariant shrinking of $S\lab$ around $t$.

For conditions (3) and (4), notice that by diagrams~\eqref{loc 3.2} and \eqref{loc 3.3} and the properties of d-critical structures, $[T\lab /G]$ is a d-critical stack and thus, up to possibly shrinking $S\lab$ $G$-invariantly around $t$, there exists an induced d-critical $G$-invariant quotient chart $(T\lab, S\lab, f\lab, \Phi\lab)$ for $\mM$ centered at $t \in T\lab$. Finally, since by construction $S\lab$, $T\lab$ are \'{e}tale slices of $V\lal, U\lal$ at $t$, we obtain that $(df\lab)=(\theta\lal^* df\lal)$ as ideals and $\theta\lal^* f\lal + I^2, f\lab + I^2$ both represent the pullback of the d-critical structure of $\mM$ to $[T\lab / G]$, which completes the proof.
\end{proof}

\begin{remark}
The above results only use the d-critical structure of $\mM$ and not the existence of its $(-1)$-shifted symplectic derived enhancement. In fact, using the property that $\mM$ is the truncation of a $(-1)$-shifted symplectic derived Artin stack $\BCal{M}$, Proposition~\ref{existence of d-crit quot} may be strengthened. Since we won't need this strengthening in this paper, we are content with just making this remark.
\end{remark}

\begin{remark}
We finally make a remark on the conditions in Definition~\ref{stack of dt type}. The first condition is there to ensure boundedness. The second condition implies the existence of quotient charts as in diagram~\eqref{strong quotient chart}, which will be necessary in order to construct the intrinsic stabilizer reduction $\tilde{\mM}$ and its good moduli space. The third condition is necessary to induce a d-critical structure on $\mM$ which is the crucial component in obtaining a semi-perfect obstruction theory for $\tilde{\mM}$.
\end{remark}

\section{Intrinsic Stabilizer Reduction} \label{sec desingularization}

In this section, we review the notions of intrinsic and Kirwan blowups and then generalize the construction of the intrinsic stabilizer reduction of a quotient stack obtained by GIT given in \cite{KLS} to the more general setting of stacks with good moduli space.

\subsection{Kirwan blowups for affine schemes} \label{kirwan blowup subsection} We recall the notion of Kirwan blowup developed in \cite{KLS}.

Suppose that $U$ is an affine scheme with an action of a reductive group $G$. For now, let us assume that $G$ is connected, as this will be the case when we take blowups throughout the paper. 

Suppose that we have an equivariant closed embedding $U \to V$ into a smooth affine $G$-scheme $V$ and let $I$ be the ideal defining $U$. Since $U \subseteq V$ is $G$-invariant, $G$ acts on $I$ and we have a decomposition $I = I^{\mathrm{fix}} \oplus I^{\mathrm{mv}}$ into the fixed part of $I$ and its complement as $G$-representations. 

Let $V^G$ be the fixed point locus of $G$ inside $V$, defined by the ideal generated by $\oO_V^{\mathrm{mv}}$, and $\pi \colon \mathrm{bl}_G(V) \to V$ the blowup of $V$ along $V^G$. 
Let $E\subseteq \mathrm{bl}_G(V)$ be its exceptional divisor and $\xi \in \Gamma( \oO_{\mathrm{bl}_G(V)}(E))$ the tautological defining equation of $E$. 

$G$-equivariance implies that (cf. \cite[Section~2.2]{KLS})
\begin{equation*}\label{xi}
\pi^{-1} (I^{\mathrm{mv}}) \subseteq \xi \cdot \oO_{\mathrm{bl}_G(V)}(-E) \subseteq \oO_{\mathrm{bl}_G(V)},
\end{equation*}
and consequently, $\xi^{-1} \pi^{-1} (I^{\mathrm{mv}}) \subseteq \oO_{\mathrm{bl}_G(V)}$, meaning that $\pi^{-1}(I^{\mathrm{mv}})$ lies in the image of the inclusion $\oO_{\mathrm{bl}_G(V)}(-E) \subseteq \oO_{\mathrm{bl}_G(V)}$, given by multiplication by $\xi$.

We define $I^{\mathrm{intr}}\subseteq\oO_{\mathrm{bl}_G(V)}$ to be 
\begin{align} \label{tilde I}
I^{\mathrm{intr}} = \text{ideal sheaf generated by } \pi^{-1}(I^{\mathrm{fix}}) \text{ and } \xi^{-1} \pi^{-1} (I^{\mathrm{mv}}).
\end{align}
\begin{definition} \emph{(Intrinsic blowup)} \label{intrinsic blow-up}
The $G$-intrinsic blowup of $U$ is the  subscheme $U^{\mathrm{intr}}\subseteq \mathrm{bl}_G(V)$ defined by the ideal $I^{\mathrm{intr}}$. 
\end{definition}

\begin{lemma}\emph{\cite[Lemma~2.6]{KLS}}\label{ideal of intrinsic}
The $G$-intrinsic blowup of $U$ is independent of the choice of $G$-equivariant embedding $U \subseteq V$,
and hence is canonical. 
\end{lemma} 

\begin{proof}
We give a brief outline of the proof and refer the reader to the proof of \cite[Lemma~2.6]{KLS} and \cite[Subsection~3.1]{KiemLi} for details. 

If $U^G = \emptyset$, then $\hU = U$ and there is nothing to show. So suppose this is not the case and $u \in U^G$.

Let $U \rightarrow V_1,\ U \rightarrow V_2$ be two equivariant embeddings of $U$ into smooth $G$-schemes. By composing with the equivariant inclusions $V_1 \times \lbrace u \rbrace \subseteq V_1 \times V_2$ and $\lbrace u \rbrace \times V_2 \subseteq V_1 \times V_2$, we may reduce to the case of a sequence of equivariant embeddings $U \rightarrow V \xrightarrow{\Phi} W$, where $V, W$ are smooth $G$-schemes.

Applying the above construction using the two embeddings of $U$, we obtain $G$-intrinsic blowups $\hU^V$, induced by the embedding $U \rightarrow V$, and $\hU^W$, induced by the composite embedding $U \rightarrow W$.

Then an identical argument as in the proof of \cite[Lemma~3.1]{KiemLi}, shows that $\Phi$ gives rise to an isomorphism $\hU^V \simeq \hU^W$, which moreover is independent of the choice of $\Phi$ interpolating between the two embeddings $U \rightarrow V$ and $U \rightarrow W$. This proves that $\hU$ is indeed canonical.
\end{proof}

Suppose $U$ is an affine $G$-scheme with an equivariant embedding into a smooth affine $G$-scheme $V$ as above. We can make sense of the notion of semistability of points in $U^{\mathrm{intr}}$ without ambiguity as follows.

We first work on the ambient scheme $V$. As it is affine, we can think of all points of $V$ as being semistable (in the usual sense of GIT). For the blowup $V^{\mathrm{intr}} := \mathrm{bl}_G (V)$, $G$ acts linearly on $E=\bP N_{V^G/V}$ with respect to the natural $G$-linearization of $\oO_E(1) := \oO_E (-E)$ and therefore we have a GIT notion of stability of points on the exceptional divisor $E$.

\begin{definition} \label{local stab} We say that $v \in V$ is \textbf{stable} if its $G$-orbit is closed in $V$ and its stabilizer finite. A point $\tilde{v} \in \mathrm{bl}_G (V)$ is \textbf{unstable} if its $G$-orbit closure meets the unstable locus of $E$. If $\tilde{v}$ is not unstable, we say that it is \textbf{semistable}.
\end{definition}

Thus for any smooth affine $G$-scheme $V$, we can define its Kirwan blowup $\hat V = \left( \mathrm{bl}_G (V) \right)^{ss}$. By \cite{Kirwan}, it satisfies $\hV^G = \emptyset$.

Now, if we have an equivariant embedding $V \to W$ between smooth $G$-schemes, then $(W^{\mathrm{intr}})^{ss} \cap V^{\mathrm{intr}} = (V^{\mathrm{intr}})^{ss}$ based on our description. Hence we may define $(U^{\mathrm{intr}})\uss := U^{\mathrm{intr}} \cap (V^{\mathrm{intr}})\uss$ for any equivariant embedding $U \to V$ into a smooth scheme $V$, which is independent of the choice of $U \to V$. 

\begin{definition} \emph{(Kirwan blowup)} \label{Kirwan blowup def}
Let $U$ be a possibly singular affine $G$-scheme with non-empty $G$-fixed locus. The Kirwan blowup of $U$ associated with $G$ is the scheme
$\hat U=(U^{\mathrm{intr}})\uss$. It satisfies $\hU^G = \emptyset$.
\end{definition}

The notion of semistability above is exactly motivated by the corresponding notion in GIT in Kirwan's original blowup procedure in \cite{Kirwan}. This can be seen by the following theorem of Reichstein, which asserts that the locus of unstable points on $\mathrm{bl}_G(V)$ is exactly the locus of unstable points in the sense of GIT for the ample line bundle $\oO_{\mathrm{bl}_G(V)}(-E)$.

\begin{theorem} \emph{\cite{Reichstein}} \label{reichstein}
Let $V$ as above and let $q \colon V \to V \git G$ denote the GIT quotient map. The unstable locus of $\mathrm{bl}_G(V)$ with respect to the line bundle $\oO_{\mathrm{bl}_G(V)}(-E)$ is the strict transform of the saturation $q^{-1}( q(V^G) )$ of $V^G$ inside $V$.
\end{theorem}

We obtain the following corollary.

\begin{corollary} \label{properties of Kirwan blowup - unstable locus, good mod} Let $\pi \colon \hat{U} \to U$ be the Kirwan blowup of an affine $G$-scheme $U$ with non-empty $G$-fixed locus. Then the quotient stack $[\hat{U}/G]$ admits a good moduli space, denoted by $q \colon [\hat{U} / G] \to \hat{U} \git G$.
\end{corollary}

\begin{proof}
Let $U \to V$ be a $G$-equivariant embedding into a smooth affine $G$-scheme $V$. By Theorem~\ref{reichstein}, $[\hV / G]$ admits a good moduli space $\hV \git G$ and therefore the same is true for the closed substack $[\hU / G]$ by Proposition~\ref{properties of good mod spaces}(6).
\end{proof}

\begin{exam} \label{Kirwan blowup example}
Suppose that $G = \bC^\ast$ is the one-dimensional torus acting on the affine plane $V = \bC^2_{x,y}$ with weights $1$ and $-1$ on the coordinates $x$ and $y$ respectively and $U \subseteq V$ is the closed $G$-invariant subscheme cut out by the ideal $I = (x^2 y, x y^2)$. 

$V^{\mathrm{intr}}$ is the blowup of $V$ along the fixed locus $V^G = \lbrace 0 \rbrace$. The unstable points are the punctured $x$-axis and $y$-axis together with the points $0, \infty$ of the exceptional divisor $\bP^1$. Thus we have that $\hV = \Spec [u ,v, v^{-1} ]$ where $G$ acts on $u, v$ with weights $1, -2$ respectively and the blowdown map $\hV \to V$ is given on coordinates by $x \mapsto u, y \mapsto uv$.

$\hU \subseteq \hV$ is the closed subscheme cut out by the ideal $(u^2)$, so that 
$$\hU = \Spec \left( \bC[u,v,v^{-1}] / (u^2) \right).$$
\end{exam}

We now explain how one can proceed if $G$ is not connected. 

Suppose that $U \to V$ is a $G$-equivariant embedding into a smooth $G$-scheme $V$. Let, as before, $I$ be the ideal of $U$ in $V$. Let $G_0$ be the connected component of the identity. This is a normal, connected subgroup of $G$ of finite index. Let $I = I^{\mathrm{fix}} \oplus I^{\mathrm{mv}}$ be the decomposition of $I$ into fixed and moving parts with respect to the action of $G_0$. Using the normality of $G_0$, we see that the fixed locus $V^{G_0}$ is a closed, smooth $G$-invariant subscheme of $V$ and also $I^{\mathrm{fix}},\ I^{\mathrm{mv}}$ are $G$-invariant. 

Let $\pi \colon \mathrm{bl}_{V^{G_0}} V \to V$ be the blowup of $V$ along $V^{G_0}$ with exceptional divisor $E$ and local defining equation $\xi$. Then, as before, take $I^{\mathrm{intr}}$ to be the ideal generated by $\pi^{-1}(I^{\mathrm{fix}})$ and $\xi^{-1}\pi^{-1}(I^{\mathrm{mv}})$. Everything is $G$-equivariant and we define $U^{\mathrm{intr}}$ as the subscheme of $\mathrm{bl}_{V^{G_0}} V$ defined by the ideal $I^{\mathrm{intr}}$.

Finally, we need to delete unstable points. By the Hilbert-Mumford criterion (cf. \cite[Theorem~2.1]{MFK}) it follows that semistability on $E$ with respect to the action of $G$ is the same as semistability with respect to the action of $G_0$, since every $1$-parameter subgroup of $G$ factors through $G_0$, and hence we may delete unstable points exactly as before and define the Kirwan blowup $\hU$.

One may check in a straightforwardly analogous way that this has the same properties (and intrinsic nature). It is obvious that if $G$ is connected we recover Definition~\ref{Kirwan blowup def}.
\medskip

\subsection{Kirwan blowups for GIT schemes} \label{Kirwan blowups GIT scheme subsection} So far, we have defined Kirwan blowups for affine $G$-schemes (with non-empty fixed locus of the $G$-action). We now generalize the construction to any $G$-scheme with a $G$-action coming from GIT (see Definition~\ref{kirwan blowup of any scheme} below) and the associated quotient stack, as performed in \cite[Subsection~2.4]{KLS}.  

We begin with an affine $G$-scheme whose $G$-fixed locus can now be possibly empty. Let $\mM = [X / G]$ be a quotient Artin stack, where $X$ is an affine $G$-scheme. 

Fix a $G$-equivariant closed embedding $X \to Y$, where $Y$ is a smooth affine $G$-scheme, and write $\yY = [ Y / G]$. Let $d$ be the maximum dimension of the stabilizers of closed points $x \in X$. By possibly equivariantly shrinking $Y$, we may assume that the maximum dimension of stabilizers of closed points of $Y$ is also equal to $d$.

Kirwan's partial desingularization, combined with Theorem~\ref{reichstein}, yields a smooth $G$-scheme $\hat{Y}$ by blowing up the locus $Y_d \subseteq Y$ consisting of points satisfying $\dim G_y = d$ and deleting unstable points. $\hat{Y}$ satisfies the crucial property that its maximum stabilizer dimension is strictly less than $d$. 

Now, for any point $x \in X$ with closed $G$-orbit and $\dim G_x = d$, let $S_x \subseteq Y$ be an \'{e}tale slice for the $G$-orbit of $x$ in $Y$. The Kirwan blowups $\hat{S}_x$ (cf. Definition~\ref{Kirwan blowup def}) of the slices $S_x$ associated with $G_x$ provide natural strongly \'{e}tale morphisms $\Psi_x \colon \hat{S}_x \times_{G_x} G \to \hat{Y}$. The joint morphism 
$$\Psi := \coprod_{x \in X_d} \Psi_x \colon \coprod_{x \in X_d} \hat{S}_x \times_{G_x} G \lr \hat{Y}$$
is an \'{e}tale cover of an open neighbourhood $\hat{Y}'$ of the restriction $E \times_{Y_d} X_d$ of the exceptional divisor $E$ of $\hat{Y}$ over the locus $X_d \subseteq Y_d$.

Let $T_x = X \cap S_x$. Then $T_x$ is an \'{e}tale slice for the $G$-orbit of $x$ in $X$ and we have the Kirwan blowup $\hat{T}_x$ using the closed $G_x$-equivariant embedding $T_x \subseteq S_x$. In \cite[Proposition~2.12]{KLS}, it is shown that the closed embedding 
$$\coprod_{x \in X_d}  \hat{T}_x \times_{G_x} G \lr \coprod_{x \in X_d}  \hat{S}_x \times_{G_x} G$$
satisfies \'{e}tale descent with respect to the morphism $\Psi$ and thus defines a closed subscheme $\hat{X}' \subseteq \hat{Y}'$. $\hat{X}'$ glues naturally with the pre-image of $X - X_d$ in $\hat{Y}$ under the blowup morphism $\hat{Y} \to Y$, defining a closed $G$-equivariant subscheme $\hat{X}$ inside $\hat{Y}$.

Lemma~\ref{ideal of intrinsic} and the above discussion on stability imply that these are canonical and independent of any choices made regarding the embedding $X \to Y$ and \'{e}tale slices $S_x$.
\medskip

If $X$ is not affine, but is the GIT semistable locus of a $G$-scheme with a $G$-linearized ample line bundle, then it admits a cover by $G$-invariant Zariski open, saturated, affine subschemes $U_1, ..., U_n \subseteq X$. Here, as in Definition~\ref{definition of saturated open substack}, saturated means that the $U_i$ fit in Cartesian squares
\begin{align*}
    \xymatrix{
    [U_i / G] \ar[r] \ar[d] & X \ar[d] \\
    U_i \git G \ar[r] & X \git G.
    }
\end{align*}

We may thus apply the Kirwan blowup to each $U_i$. The Kirwan blowups $\hat{U}_i$ glue to the Kirwan blowup $\hat{X}$ of $X$ and we are led to the following definition. 

\begin{definition} \label{kirwan blowup of any scheme}
Suppose that $[X/G]$ is a GIT quotient stack. Then we call $\hat{X}$ and $\hat{\mM} = [\hat{X} / G]$ the Kirwan blowups of $X$ and $\mM = [X / G]$ respectively. 
\end{definition}

\begin{remark}
When $X$ is not affine, we could also argue globally as in the affine case using an equivariant closed embedding $X \to Y$ into a smooth $G$-scheme $Y$, afforded by GIT. However, the more local argument using a cover by saturated, affine, $G$-invariant open subschemes is closer in spirit to the proof of Theorem-Construction~\ref{intrinsic stab red of stack with gms}.
\end{remark}

\subsection{Intrinsic stabilizer reduction for GIT quotient stacks}

We quickly recall how the construction of the intrinsic stabilizer reduction works for GIT quotient stacks and some useful properties in that setting.

Let $X$ be a $G$-scheme coming from GIT and $\mM = [ X / G]$ the associated quotient stack.

Write $\mM_0 = \mM$ and $X_0 = X$. Taking Kirwan blowups, we obtain $X_1 = \hat{X}$ and $\mM_1 = [\hat{X} / G]$. The $G$-action on $\hat{X}$ is still $G$-linearized, so $\mM_1$ is a GIT quotient stack of lower maximum stabilizer dimension compared to $\mM_0$. Thus we may iterate the procedure to produce a canonical sequence of GIT quotient stacks
$$\mM_0 = [X / G],\ \mM_1 = [X_1 / G],\ \ldots,\ \mM_\ell = [X_\ell / G],$$
where $\mM_\ell$ is a Deligne--Mumford stack, since at each step the maximum stabilizer dimension strictly decreases.

\begin{definition}
$\tilde{\mM} := \mM_\ell$ is called the intrinsic stabilizer reduction of $\mM$.
\end{definition}

\subsection{Intrinsic stabilizer reduction for stacks with good moduli spaces}

We now generalize the procedure of intrinsic stabilizer reduction to the case of stacks with good moduli spaces.

The following properties of Kirwan blowups will be useful, so we record them separately in a proposition before our main construction.

\begin{proposition} \label{prop 4.9}
Let $\varphi \colon G_1 \to G_2$ be a surjective homomorphism between reductive groups, and $X_i$, $i=1,2$, be affine $G_i$-schemes. 

Suppose that $f \colon X_1 \to X_2$ is a smooth morphism which is equivariant with respect to $\varphi$ and such that the induced morphism $\underline{f} \colon [X_1 / G_1 ] \to [X_2 / G_2]$ on quotient stacks is strongly \'{e}tale (cf. Definition~\ref{definition of strongly etale}).

Then there is a canonically induced smooth, affine morphism $\hat{f} \colon \hat{X}_1 \to \hat{X}_2$, equivariant with respect to $\varphi$, such that the associated quotient morphism $\underline{\hat{f}} \colon [\hat{X}_1 / G_1] \to [\hat{X}_2 / G_2]$ is strongly \'{e}tale and fits in a Cartesian square
\begin{align*}
    \xymatrix{
    [\hat{X}_1 / G_1] \ar[d] \ar[r]^-{\underline{\hat{f}}} & [\hat{X}_2 / G_2] \ar[d] \\
    [X_1 / G_1] \ar[r]_-{\underline{f}} & [X_2 / G_2].}
\end{align*}
\end{proposition}

\begin{proof}
Let $x_1 \in X_1$ a point with closed $G_1$-orbit and of maximum stabilizer dimension and $x_2 = f(x_1) \in X_2$. Since $\underline{f}$ is strongly \'{e}tale, the $G_2$-orbit of $x_2$ is also closed, by Proposition~\ref{properties of good mod spaces}(5), and $x_2$ is of maximum stabilizer dimension.

Up to $G_1$-equivariant shrinking around $x_1$, we can find smooth, affine $G_i$-schemes $Y_i$, $i=1,2$, and a smooth morphism $g \colon Y_1 \to Y_2$, equivariant with respect to $\varphi$, fitting in a commutative diagram
\begin{align} \label{loc 4.2}
    \xymatrix{
    X_1 \ar[r] \ar[d]_-{f} & Y_1 \ar[d]^-{g} \\
    X_2 \ar[r] & Y_2
    }
\end{align}
where the horizontal arrows are equivariant closed embeddings. We may additionally assume that the induced morphism $\underline{g} \colon [Y_1 / G_1] \to [Y_2 / G_2]$ is strongly \'{e}tale. This is because we can initially take $g$ to map $x_1$ to $x_2$, be stabilizer-preserving at $x_1$ and $\underline{g}$ to be \'{e}tale at $x_1$, so using the Fundamental Lemma \cite[Theorem~6.10]{AlperFundLemma} at $x_1$, after possible further shrinking around $x_1$, we may take $\underline{g}$ to be strongly \'{e}tale.

By the intrinsic nature of intrinsic blowups, using an identical argument as the one outlined in the proof of Lemma~\ref{ideal of intrinsic}, one may consider diagram~\eqref{loc 4.2} as $x_1$ varies to obtain a canonical morphism $f^{\mathrm{intr}} \colon X_1^{\mathrm{intr}} \to X_2^{\mathrm{intr}}$.

By construction, $f^{\mathrm{intr}}$ is stabilizer-preserving and maps closed $G_1$-orbits to closed $G_2$-orbits, as this is the case for $g^{\mathrm{intr}}$. It is moreover smooth and the associated quotient morphism $\underline{f}^{\mathrm{intr}} \colon [X_1^{\mathrm{intr}} / G_1] \to [X_2^{\mathrm{intr}} / G_2]$ is \'{e}tale. This can be seen as follows: Let $T$ be an affine, \'{e}tale slice for the $G_1$-orbit of $x_1$ in $X_1$, whose stabilizer we denote by $H$. By the assumptions on the morphism $f$, we have a commutative diagram
\begin{align*}
    \xymatrix{
    T \times_H G_1 \ar[d]_{\id \times \phi} \ar[r] & X_1 \ar[d]^-{f} \\
    T \times_H G_2 \ar[r] & X_2
    }
\end{align*}
and the corresponding commutative diagram at the level of quotient stacks
\begin{align*}
    \xymatrix{
    [T/H] \ar[d]_{\id } \ar[r] & [X_1 / G_1] \ar[d]^-{\underline{f}} \\
    [T/H] \ar[r] & [X_2 / G_2]
    }
\end{align*}
has all \'{e}tale arrows.

By definition, these give rise to identical diagrams with $T, X_1, X_2$ replaced by their intrinsic blowups. Varying $x_1$, this immediately implies the smoothness and affineness of ${f^{\mathrm{intr}}}$ and \'{e}taleness of $\underline{f}^{\mathrm{intr}}$. Moreover, a similar, straightforward local computation, using the definition of intrinsic blowup and the arguments of the proof of Lemma~\ref{ideal of intrinsic} and \cite[Lemma~3.1]{KiemLi}, shows that the natural map $X_1^{\mathrm{intr}} \to X_1 \times_{X_2} X_2^{\mathrm{intr}}$ is an isomorphism.

Now Corollary~\ref{properties of Kirwan blowup - unstable locus, good mod} lets us pass from intrinsic blowups to Kirwan blowups.

Finally, since $\underline{\hat{f}}$ is \'{e}tale, stabilizer-preserving and maps closed points to closed points, applying the Fundamental Lemma~\cite[Theorem~6.10]{AlperFundLemma} at all closed points of $[\hat{X}_1/G_1]$, we obtain that it is strongly \'{e}tale. The existence of the isomorphism $[\hat{X}_1 / G_1] \to [X_1 / G_1] \times_{[X_2 / G_2]} [\hat{X}_2 / G_2]$ is a direct consequence of the isomorphism $X_1^{\mathrm{intr}} \to X_1 \times_{X_2} X_2^{\mathrm{intr}}$ and the definition of stability (Definition~\ref{local stab}). This concludes the proof.
\end{proof}

Before extending the notion of Kirwan blowup to stacks with good moduli space, it will be useful to generalize our definition of stability to this setting (cf. Definition~\ref{local stab}). 

\begin{definition}
Let $\mM$ be an Artin stack of finite type and $q \colon \mM \to M$ a good moduli space morphism of finite type with affine diagonal. The prestable locus $\mM^{ps} \subseteq \mM$ of $\mM$ is the (saturated) open substack whose closed points are the closed points $x \in \mM$ satisfying $ |q^{-1}(q(x))| = \lbrace x \rbrace $. The stable locus $\mM^s \subseteq \mM$ is the saturated, open substack consisting of prestable points with finite stabilizer group.
\end{definition}

This agrees with Definition~\ref{local stab} by the following proposition, which also explains why the stable locus is open and saturated, since this is a property that can be checked strongly \'{e}tale locally.

\begin{proposition}
With the same notation as above, for any strongly \'{e}tale morphism $[U/G] \xrightarrow{\Phi} \mM$, we have that $\Phi^{-1} (\mM^s) = [U^s / G]$, where $U^s$ is the open subscheme of $u \in U$ with closed $G$-orbit and finite stabilizer. 
\end{proposition}

\begin{proof}
Any strongly \'{e}tale morphism
$$
\xymatrix{
[U/G] \ar[r]^-\Phi \ar[d]_-{q_U} & \mM \ar[d]^-{q} \\
U \git G \ar[r] & M
}
$$
induces an isomorphism of fibers $q_U^{-1}(q_U(u)) \simeq q^{-1}(q(x))$ whenever $\Phi(u) = x$ and we identify $u \in U$ with its $G$-orbit as a point $u \in [U/G]$. The conclusion is immediate by Proposition~\ref{properties of good mod spaces}(5) and the fact that $\Phi$ is stabilizer-preserving.
\end{proof}

We can now construct the Kirwan blowup of an Artin stack with good moduli space.

\begin{thm-constr} \label{Kirwan blowup of stack with GMS} Let $\mM$ be an Artin stack of finite type over $\bC$ with affine diagonal, which is not Deligne--Mumford. Moreover, suppose that $q \colon \mM \to M$ is a good moduli space morphism with $q$ of finite type and with affine diagonal. 

Then there exists a canonical Artin stack $\hat{\mM}$, called the Kirwan blowup of $\mM$, together with a morphism $\pi \colon \hat{\mM} \to \mM$, such that: 
\begin{enumerate}
    \item $\hat{\mM}$ is of finite type over $\bC$, has affine diagonal and admits a good moduli space morphism $\hat{q} \colon \hat{\mM} \to \hat{M}$ with affine diagonal.
    \item The maximum stabilizer dimension of closed points in $\hat{\mM}$ is strictly smaller than that of $\mM$.
    \item For any strongly \'{e}tale morphism $[U/G] \to \mM$, the base change $\hat{\mM} \times_{\mM} [U/G]$ is naturally isomorphic to the Kirwan blowup $[\hU / G]$. 
    \item $\pi |_{\pi^{-1}(\mM^s)}$ is an isomorphism.
\end{enumerate} 

As is the case with Kirwan blowups of schemes, $\hat{\mM}$ is a semistable locus $\hat{\mM} = (\mM^{\mathrm{intr}})^{ss}$, an open substack of a canonical Artin stack $\pi^{\mathrm{intr}} \colon \mM^{\mathrm{intr}} \to \mM$, called the intrinsic blowup of $\mM$.
\end{thm-constr}
\begin{proof}
Let $\mM^{max}$ be the substack of $\mM$ whose points have stabilizers of the maximum possible dimension. This is a closed substack of $\mM$ with good moduli space $M^{max}$ (cf. \cite[Appendix B]{EdidinRydh} for details on the algebraic structure of $\mM^{max}$). For any closed point $x \in \mM^{max}$, applying Theorem-Definition~\ref{quotient chart existence thm-def}, we have a Cartesian diagram
\begin{align}
\xymatrix{
\left[ U_x / G_x \right] \ar[r]^-{\Phi_x} \ar[d]_{q_x} & \mM \ar[d]^-q \\
U_x \git G_x \ar[r] & M.
}
\end{align}
The morphisms $\Phi_x$ are affine, strongly \'{e}tale and cover the locus $\mM^{max}$. We may take the Kirwan blowup of each quotient stack $[U_x / G_x]$ to obtain good moduli space morphisms $[\hU_x / G_x] \to \hU_x \git G_x$.

We need to check that these glue to give a stack $\hat{\mM}$ with a universally closed projection $\hat{\mM} \to \mM$ and a good moduli space $\hat{\mM} \to \hat{M}$ satisfying the same conditions as $\mM$ and its good moduli space morphism. By the properties of the Kirwan blowup, the maximum stabilizer dimension of $\hat{\mM}$ will be lower than that of $\mM$.

Suppose $x, y$ are two closed points of $\mM$ such that $G_x, G_y$ are of maximum dimension. We obtain a Cartesian diagram of stacks
\begin{align} \label{square1}
\xymatrix{
[U_x \times_\mM U_y / (G_x \times G_y)] \ar[r] \ar[d] & [U_y / G_y] \ar[d] \\
[U_x / G_x] \ar[r] & \mM
} 
\end{align}
where $U_{xy} := U_x \times_\mM U_y$ is an affine scheme. This is due to the Cartesian diagram 
\begin{align*}
\xymatrix{
U_x \times_{\mM} U_y \ar[r] \ar[d] & \mM \ar[d]^-{\Delta_{\mM}} \\
U_x \times U_y \ar[r] & \mM \times \mM
}
\end{align*}
and the fact that $\mM$ has affine diagonal.

Using Proposition~\ref{prop 4.9}, we obtain a diagram
\begin{align} \label{fiber prod}
\xymatrix{
& [\hat{U}_{xy} / (G_x \times G_y)] \ar[dl] \ar[dr]\\
[\hU_x / G_x] & & [\hU_y / G_y]
} 
\end{align}
with affine, strongly \'{e}tale arrows and we have, moreover, canonical associated isomorphisms $[\hat{U}_{xy} / (G_x \times G_y)] \to [\hU_x / G_x] \times_\mM [U_y / G_y]$ and $[\hat{U}_{xy} / (G_x \times G_y)] \to [U_x / G_x] \times_\mM [\hU_y / G_y]$. 

Using the charts $[\hU_x / G_x]$ together with a cover of $\mM \setminus q^{-1}(M^{max})$, we therefore obtain an atlas for a stack $\hat{\mM}$ with a map to $\mM$. By the canonical isomorphisms of the previous paragraph, $\hat{\mM}$ is independent of the particular choices of charts for $\mM$.
 
Since the arrows in diagram~\eqref{fiber prod} are strongly \'{e}tale, we obtain a corresponding diagram of \'{e}tale arrows at the level of good moduli spaces of the Kirwan blowups
\begin{align*}
\xymatrix{
& [\hat{U}_{xy} / (G_x \times G_y)] \ar[dl] \ar[dr] \ar[dd]\\
[\hU_x / G_x] \ar[dd] & & [\hU_y / G_y] \ar[dd] \\
& \hat{U}_{xy} \git (G_x \times G_y) \ar[dl] \ar[dr]\\
\hU_x \git G_x & & \hU_y \git G_y
} 
\end{align*}
where both rhombi are Cartesian. 

By Proposition~\ref{properties of good mod spaces}, $\mM \setminus q^{-1} (M^{max}) \to M \setminus M^{max}$ is a good moduli space morphism. Hence the morphisms $[\hU_x / G_x] \to \hU_x \git G_x$ for all $x \in \mM^{max}$ together with an atlas of $\mM \setminus q^{-1} (M^{max})$ glue to give a morphism $\hat{\mM} \to \hat{M}$. By Proposition~\ref{properties of good mod spaces} again, this is a good moduli space morphism.

$\hat{\mM}$ has affine diagonal since we have a cartesian diagram 
$$\xymatrix{
[\hU_{xy} / (G_x \times G_y)] \ar[d] \ar[r] & \hat{\mM} \ar[d] \\
[\hU_x / G_x] \times [\hU_y / G_y] \ar[r] & \hat{\mM} \times \hat{\mM}
}$$
where the lower horizontal arrows give an \'{e}tale cover of $\hat{\mM} \times \hat{\mM}$ and the left vertical arrow is affine.

To see that $\hat{\mM} \to \hat{M}$ also has affine diagonal, we consider the diagram
\begin{align*}
\xymatrix{
\hat{\mM} \ar[r] \ar[dr] & \hat{\mM} \times_{\hat{M}} \hat{\mM} \ar[r] \ar[d] & \hat{M} \ar[d] \\
 & \hat{\mM} \times \hat{\mM} \ar[r] & \hat{M} \times \hat{M}
}
\end{align*}
where the right square is cartesian. The diagonal of $\hat{M}$ is separated \cite[\href{https://stacks.math.columbia.edu/tag/04YQ}{Tag 04YQ}]{stacks-project}. Therefore, since the diagonal of $\hat{\mM}$ is affine, it follows by the usual cancellation property that $\hat{\mM} \to \hat{M}$ has affine diagonal.

$\hat{\mM}$, $\hat{M}$ and the morphism $\hat{\mM} \to \hat{M}$ thus have the same properties as $\mM$, $M$ and $\mM \to M$, as desired. This establishes property (1) of the statement.

Property (3) is true by construction and properties (2) and (4) can be checked strongly \'{e}tale locally. But by the definition of Kirwan blowups of quotient stacks, they are both true in that case.

Finally, the last assertion about the intrinsic blowup $\mM^{\mathrm{intr}}$ can be shown using the same arguments. Since we won't be needing $\mM^{\mathrm{intr}}$, we leave the details to the reader.
\end{proof}

\begin{remark}
The fact that $\mM$ and $q \colon \mM \to M$ have affine diagonal is not crucial for the above proof to go through. In fact, by \cite[Theorem~13.1]{AHR2} and \cite[Corollary~13.11]{AHR2}, it is enough to assume that $\mM$ has affine stabilizers and separated diagonal. Since the stacks we are interested in will have affine diagonal, we make this assumption for convenience of presentation. 
\end{remark}

As before, repeatedly applying the operation of Kirwan blowup we obtain the instrisic stabilizer reduction. By construction of the Kirwan blowup, the stable locus is preserved at each step and hence stays invariant under the whole procedure.

\begin{theorem} \label{intrinsic stab red of stack with gms}
Let $\mM$ be an Artin stack of finite type over $\bC$ with affine diagonal, which is not Deligne--Mumford. Moreover, suppose that $q \colon \mM \to M$ is a good moduli space morphism with $q$ of finite type and with affine diagonal. Then there exists a canonical Deligne--Mumford stack $\tilde{\mM}$, called the intrinsic stabilizer reduction of $\mM$, together with a morphism $\pi \colon \tilde{\mM} \to \mM$ which restricts to an isomorphism over the stable locus $\mM^s \subseteq \mM$. Moreover, $\tilde{\mM}$ admits a good moduli space $\tilde{M}$ and the good moduli space morphism $\tilde{q} \colon \tilde{\mM} \to \tilde{M}$ is proper.
\end{theorem}

\begin{proof} We perform a sequence of Kirwan blowups by applying Theorem-Costruction~\ref{Kirwan blowup of stack with GMS} iteratively to get
$$\mM_0 = \mM,\ \mM_1 = \hat{\mM}_1,\ \ldots,\ \mM_\ell = \hat{\mM}_{\ell-1}.$$
At each step, the maximum stabilizer dimension strictly decreases, so we terminate as soon as $\mM_{\ell}$ has finite stabilizers and is thus Deligne--Mumford. We define $\tilde{\mM} := \mM_\ell$.

$\tilde{\mM} \to \tilde{M}$ admits a strongly \'{e}tale cover by morphisms of the form $[U/G] \to U \git G$ where $G$ is finite, so its diagonal is finite (cf. \cite[Proposition~0.8]{MFK}). Thus $\tilde{\mM} \to \tilde{M}$ is separated and by Proposition~\ref{properties of good mod spaces} also universally closed, hence proper.
\end{proof}

\begin{remark}
If $\mM^s$ is non-empty, then we are guaranteed that the intrinsic stabilizer reduction $\tilde{\mM}$ is non-empty. While this is sufficient, it is not necessary, as can be seen for the quotient stack $\mM = [U / G]$ in Example~\ref{Kirwan blowup example}. 

In general, if at some step of the procedure the stack $\mM_i$ is a gerbe over a Deligne--Mumford stack (i.e. $\mM_i$ has constant stabilizer dimension and every point is prestable), then $\hat{\mM}_i$ is empty. We could elect to terminate the procedure at this step instead and define $\tilde{\mM} := \mM_i$. We choose not to do so as for our purpose of defining Donaldson--Thomas invariants, setting $\tilde{\mM} = \emptyset$ in this case seems more geometrically motivated and better aligned with the BPS invariants considered in \cite{JoyceSong} and \cite{DavMein}.
\end{remark}

\begin{remark}
Edidin-Rydh have also developed a blowup procedure for stacks with good moduli spaces in \cite{EdidinRydh}. For smooth stacks, our stabilizer reduction is the same as theirs. For singular stacks, Kirwan blowups can be phrased in their language of saturated blowups, however the resolution they obtain is a closed substack of the one here. Nevertheless, our construction is closely related to the Edidin-Rydh resolution of stabilizers performed at the level of derived stacks, using the recently developed notion of derived blowups in \cite{Hekking} and the results of \cite{HRS}.
\end{remark}

\section{Obstruction Theory} \label{sec obstruction theory}

In this section, we first recall the basic principles of semi-perfect obstruction theories and then explain how to construct such a gadget on the intrinsic stabilizer reduction of a stack of DT type. Our discussion goes through the notion of a local model and its obstruction theory, as used in \cite{KLS}.

\subsection{Semi-perfect obstruction theory} 
This subsection contains background material about semi-perfect obstruction theories and their induced virtual cycles, as developed in \cite{LiChang}.

Let $U \to C$ be a morphism, where $U$ is a scheme of finite type and $C$ denotes a smooth quasi-projective scheme, which will typically be either a point or a smooth quasi-projective curve. We first recall the definition of 
perfect obstruction theory \cite{BehFan, LiTian}.

\begin{definition} \emph{(Perfect obstruction theory \cite{BehFan})} \label{Perf obs th}
A (truncated) perfect (relative) obstruction theory consists of a morphism $\phi \colon E \to L_{U/C}^{\geq -1}$ in $D^b(\Coh U)$ such that
\begin{enumerate}
\item $E$ is of perfect amplitude, contained in $[-1,0]$.
\item $h^0(\phi)$ is an isomorphism and $h^{-1}(\phi)$ is surjective.
\end{enumerate}
We refer to $\Ob_\phi := \Hone(E^\vee)$ as the obstruction sheaf of $\phi$.
\end{definition}

\begin{definition} \emph{(Infinitesimal lifting problem)} \label{Inf lift prob}
Let $\iota \colon \Delta \to \bar{\Delta}$ be an embedding with $\bar{\Delta}$ local Artinian, such that $I \cdot \fm = 0$ where $I$ is the ideal of $\Delta$ and $\fm$ the closed point of $\bar{\Delta}$. We call $(\Delta, \bar{\Delta}, \iota, \fm)$ a small extension. Given a commutative square
\begin{align}
\xymatrix{
\Delta \ar[r]^g \ar[d]^\iota & U \ar[d]\\
\bar{\Delta} \ar[r] \ar@{-->}[ur]_{\bar{g}} & C
}
\end{align}
such that the image of $g$ contains a point $p \in U$, the problem of finding $\bar{g} \colon \bar{\Delta} \to U$ making the diagram commutative is the ``infinitesimal lifting problem of $U/C$ at $p$".
\end{definition}

\begin{definition} \emph{(Obstruction space)} \label{Obs spaces}
For a point $p \in U$, the intrinsic obstruction space to deforming $p$ is $T_{p, U/ C}^1 := \Hone \left( (L_{U/C}^{\geq -1})^\vee \vert_p \right)$. The obstruction space with respect to a perfect obstruction theory $\phi$ is $\OB(\phi,p) := \Hone( E^\vee \vert_p )$.
\end{definition}

Given an infinitesimal lifting problem of $U/C$ at a point $p$, there exists by the standard theory of the cotangent complex a canonical element 
\begin{align}
\omega \left( g, \Delta, \bar{\Delta} \right) \in \Ext^1 \left( g^{\ast} L_{U/C}^{\geq -1} \vert_p, I\right) = T^1_{p, U/C} \otimes_\bC I
\end{align} 
whose vanishing is necessary and sufficient for the lift $\bar{g}$ to exist. 

\begin{definition} \emph{(Obstruction assignment)} \label{Obs assignment}
For an infinitesimal lifting problem of $U / C$ at $p$ and a perfect obstruction theory $\phi$ the obstruction assignment at $p$ is the element
\begin{align}
ob_U(\phi,g,\Delta,\bar{\Delta}) = h^1(\phi^\vee) \left( \omega \left( g, \Delta, \bar{\Delta} \right) \right) \in \OB(\phi,p) \otimes_\bC  I.
\end{align}
\end{definition}

\begin{definition} \label{Same obs assign}
Let $\phi \colon E \to L_{U / C}^{\geq -1}$ and $\phi' \colon E' \to L_{U / C}^{\geq -1}$ be two perfect obstruction theories and $\psi \colon \Ob_\phi \to \Ob_{\phi'}$ be an isomorphism. We say that the obstruction theories give the same obstruction assignment via $\psi$ if for any infinitesimal lifting problem of $U/C$ at $p$
\begin{align}
\psi \left( ob_U(\phi,g,\Delta,\bar{\Delta}) \right) = ob_U(\phi',g,\Delta,\bar{\Delta}) \in \OB(\phi',p) \otimes_\bC I.
\end{align}
\end{definition}
We are now ready to give the definition of a semi-perfect obstruction theory.

\begin{definition} \emph{(Semi-perfect obstruction theory \cite[Definition~3.1]{LiChang})} \label{spot def}
Let $\mM\to C$ be a morphism, where $\mM$ is a DM stack, proper over $C$, of finite presentation and $C$ is a smooth quasi-projective scheme. A semi-perfect obstruction theory $\phi$ consists of an \'{e}tale covering $\lbrace U_\alpha \rbrace_{\alpha \in A}$ 
of $\mM$ and perfect obstruction theories $\phi_\alpha \colon E_\alpha \to L_{U_\alpha / C}^{\geq -1}$ such that
\begin{enumerate}
\item For each pair of indices $\alpha, \beta$, there exists an isomorphism \begin{align*}
\psi_{\alpha \beta} \colon \Ob_{\phi_\alpha} \vert_{U_{\alpha\beta}} \lra \Ob_{\phi_\beta} \vert_{U_{\alpha\beta}}
\end{align*}
so that the collection $\lbrace \Ob_{\phi\lalp}, \psi\lab \rbrace$ gives descent data of a coherent sheaf on $\mM$.
\item For each pair of indices $\alpha, \beta$, the obstruction theories $E_\alpha \vert_{U_{\alpha \beta}}$ and $E_\beta \vert_{U_{\alpha \beta}}$ give the same obstruction assignment via $\psi_{\alpha \beta}$ (as in Definition \ref{Same obs assign}).
\end{enumerate}
\end{definition}

\begin{remark}
The obstruction sheaves $\lbrace \Ob_{\phi\lalp} \rbrace_{\alpha \in A}$ glue to define a sheaf $\Ob_{\phi}$ on $\mM$. This is the obstruction sheaf of the semi-perfect obstruction theory $\phi$.
\end{remark}

Suppose now that $\mM \to C$ is as above and admits a semi-perfect obstruction theory. Then, for each $\alpha \in A$, we have
\begin{align*}
\cC_{U\lalp/C} \subset N_{U\lalp/C} = h^1/h^0 ( (L_{U\lalp/C}^{\geq -1})^\vee ) \lhook\joinrel\xrightarrow{h^1/h^0(\phi\lalp^\vee)} h^1/h^0 ( E\lalp^\vee ) \lr h^1(E\lalp^\vee),
\end{align*}
where $\cC_{U\lalp/C}$ and $N_{U\lalp/C}$ denote the intrinsic normal cone stack and intrinsic normal sheaf stack respectively, where by abuse of notation we identify a sheaf $\fF$ on $\mM$ with its sheaf stack.

We therefore obtain a cycle class $[\fc_{\phi \lalp}] \in Z_* \Ob_{\phi\lalp}$ by taking the pushforward of the cycle $[\cC_{U\lalp/C}] \in Z_* N_{U\lalp/C}$. 

\begin{thm-defi} \emph{\cite[Definition-Theorem~3.7]{LiChang}} \label{spot virtual cycle def}
Let $\mM$ be a DM stack, proper over $C$, of finite presentation and $C$ a point or a smooth quasi-projective curve, such that $\mM \to C$ admits a semi-perfect obstruction theory $\phi$. The classes $[\fc_{\phi \lalp}] \in Z_* \Ob_{\phi\lalp}$ glue to define an intrinsic normal cone cycle $[\fc_\phi] \in Z_* \Ob_\phi$. Let $s$ be the zero section of the sheaf stack $\Ob_\phi$. The virtual cycle of $\mM$ is defined to be
\begin{align*}
[\mM, \phi]^{\mathrm{vir}} := s^{!} [\fc_\phi] \in A_* \mM,
\end{align*}
where $s^{!} \colon Z_* \Ob_\phi \to A_* \mM$ is the Gysin map. This virtual cycle satisfies all the usual properties, such as deformation invariance.
\end{thm-defi}

\begin{remark} \label{spot schemes to dm stacks}
Observe that in Definition~\ref{spot def} it is not strictly necessary to take the $U\lalp$ to be schemes. It is straightforward to generalize the definition and the construction of the associated virtual fundamental cycle to include \'{e}tale covers by Deligne--Mumford stacks.
\end{remark}

\subsection{Local models, standard forms and their blowups} \label{local model}
Let $V$ be a smooth affine $G$-scheme. The action of $G$ on $V$ induces a morphism $\fg \otimes \oO_V \rightarrow T_V$ and its dual $\sigma_V : \Omega_V \rightarrow \fg^{\vee} \otimes \oO_V$. 

\begin{setup-def} \cite[Setup-Definition~5.1]{KLS} \label{ind hyp}
Consider the quadruple $(V, F_V, \omega_V, D_V)$, where $F_V$ is a $G$-equivariant vector bundle on $V$, $\omega_V$ a $G$-invariant section with scheme-theoretic zero locus $U = \lbrace \omega_V = 0 \rbrace \sub V$ and $D_V \sub V$ an effective invariant divisor, satisfying:
\begin{enumerate}
\item $\sigma_V(-D_V)  : \Omega_V ( -D_V ) \to \fg^\vee (-D_V)$ factors through a morphism $\phi_V$ as shown
\begin{align} \label{5.23}
\Omega_V(-D_V) \lra F_V \xrightarrow{\phi_V} \fg^\vee(-D_V);
\end{align}
\item The composition $\phi_V \circ \omega_V$ vanishes identically;
\item Let $R$ be the identity component of the stabilizer group of a closed point in $V$ with closed orbit. Let $V^R$ denote the fixed point locus of $R$. Then $\phi_V |_{V^R}$ composed with the projection $\fg^\vee(-D_V) \rightarrow \fr^\vee(-D_V)$ is zero, where $\fr$ is the Lie algebra of $R$.
\end{enumerate}
We say that the data $$\Lambda_U = (U, V, G, F_V, \omega_V, D_V, \phi_V)$$ give a local model structure for $U$. Thinking of the quotient stack, we say that $[U/G]$ also has a local model structure and denote the data by $\Lambda_{[U/G]}$.
\end{setup-def}

\begin{remark} \label{d-crit is local model}
Note that if $f \colon V \to \bA^1$ is a $G$-invariant function on $V$, then $(U, V, G, \Omega_V, df, 0, \sigma_V)$ give a local model for $U$, being equivalent to an invariant d-critical chart $(U,V,f,i)$ for $U$. Therefore, an invariant d-critical locus is a particular case of a local model.
\end{remark}

Now, let $\Lambda_U= (U, V, G, F_V, \omega_V, D_V, \phi_V)$ define a local model structure on $U$. Since $G$ is reductive, we have a splitting
\begin{align*}
F_V |_{V^G} = F_V |_{V^G}^{\mathrm{fix}} \oplus F_V |_{V^G}^{\mathrm{mv}}.
\end{align*}
Let now $\pi \colon \hV \to V$ be the Kirwan blowup of $V$. Define $F_{\hV}$ as the kernel of the composite morphism $\pi^* F_V \to \pi^* \left( F_V |_{V^G} \right) \to \pi^* \left( F_V |_{V^G}^{\mathrm{mv}} \right)$ so that we have an exact sequence
$$0 \lra F_{\hV} \lra \pi^* F_V \lra \pi^* \left( F_V |_{V^G}^{\mathrm{mv}} \right) \lra 0.$$
By equivariance, $\pi^* \omega_V$ maps to zero under the second map and hence induces an invariant section $\omega_{\hV}$ of $F_{\hV}$. A local computation shows the following.

\begin{proposition}
The zero locus of $\omega_{\hV}$ is the Kirwan blowup $\hU$ of $U$ with respect to $G$. 
\end{proposition}

\begin{proof}
This is Proposition~2.11 in \cite{KLS}.
\end{proof}

Local model structures are preserved under the operations of Kirwan blowups and taking \'{e}tale slices.

\begin{proposition} \label{propagation of local model structure}
Let $(U, V, G, F_V, \omega_V, D_V, \phi_V)$ be the data of a local model. Then there exist induced data $(\hU, \hV, G, F_{\hV}, \omega_{\hV}, D_{\hV}, \phi_{\hV})$ of a local model on the Kirwan blowup $\hU$. The same is true for any \'{e}tale slice $T$ of a closed point of $\hU$ with closed $G$-orbit.
\end{proposition}

\begin{proof}
This is \cite[Lemma~5.3]{KLS}. The fact that $F_{\hV}$ is locally free is discussed shortly after Definition~2.10 in \cite{KLS}. Denoting the Kirwan blowup map $\hV \to V$ by $\pi$, the divisor $D_{\hV}$ is equal to $\pi^* D_V + 2E$, where $E$ is the exceptional divisor of $\hV$. Finally, $\phi_{\hV}$ is induced by $\phi_V$ through a diagram chase and the definition of $F_{\hV}$. The reasoning for \'{e}tale slices is similar.
\end{proof}

\subsection{Obstruction theory of local model} Let $$\Lambda_U = (U,V,G, F_V,\omega_V, D_V, \phi_V)$$ be the data of a local model structure. Consider the sequence
\begin{align}
\kK_U := [ \fg \to T_V |_U \xrightarrow{\left( d_V \omega_V^\vee \right)^\vee} F_V |_U \stackrel{\phi_V}{\lr} \fg^\vee(-D_V) ].
\end{align}
By \cite[Lemma~5.5]{KLS}, this defines a $G$-equivariant perfect complex on $U$ that descends to a perfect complex on $[U/G]$, which we denote by $\kK_{[U/G]}$. 

We have a morphism 
\begin{align} \label{loc 5.7}
    E_{[U/G]} := \kK_{[U/G]}^\vee \lr L_{[U/G]}^{\geq -1} 
\end{align}
given by the diagram
\begin{align} \label{loc 5.8}
    \xymatrix{
    \fg(D_V) \ar[r]^-{\phi_V^\vee} \ar[d] & F_V^\vee|_U \ar[r]^-{d_V \omega_V^\vee} \ar[d]_-{\omega_V^\vee} & \Omega_V|_U \ar[r] \ar[d] & \fg^\vee \ar[d] \\
    0 \ar[r] & I/I^2 \ar[r]_{d_V} & \Omega_V|_U \ar[r] & \fg^\vee,
    }
\end{align}
where $I$ denotes the ideal sheaf of $U$ in $V$. Observe that the diagram commutes, since by Setup-Definition~\ref{ind hyp}(2), the composition $\phi_V \circ \omega_V$ is identically zero.

Let $U^s, V^s$ denote the stable loci of $U$ and $V$ respectively (cf. Definition~\ref{local stab}). We thus have that $\sigma_V|_{V^s} \colon \Omega_V|_{V^s} \to \fg^\vee$ is surjective. By Setup-Definition~\ref{ind hyp}, it follows that $\phi_V|_{U^s}$ is surjective and hence its dual is injective with locally free cokernel. 

Therefore, $\kK_{[U^s/G]} := \kK_{[U/G]} |_{[U^s / G]}$ is a two-term perfect complex, whose dual $E_{[U^s / G]}$ gives a perfect obstruction theory on the DM stack $[U^s / G]$ via the diagrams~\eqref{loc 5.7} and \eqref{loc 5.8}.

\begin{definition} \label{pot on stable locus}
Let $(U,V,G,F_V,\omega_V, D_V, \phi_V)$ be the data of a local model structure. Then $E_{[U^s / G]}$ is the induced perfect obstruction theory on the stable locus $[U^s / G]$.
\end{definition}

\subsection{Semi-perfect obstruction theory of the intrinsic stabilizer reduction of a stack of DT type} 

Let $\mM$ be an Artin stack of DT type. We now explain how to construct a semi-perfect obstruction theory on its intrinsic stabilizer reduction $\tilde{\mM}$, one of the main results of this paper.

We begin by recording and abstracting some of the data associated to a stack of DT type.

\begin{setup} \label{VC package - setup of stack of DT type}
For any closed point $x \in \mM$, by Proposition~\ref{existence of d-crit quot}, there exists a strongly \'{e}tale d-critical quotient chart for $\mM$ centered at $x$. In particular, we obtain a strongly \'{e}tale cover
\begin{align} \label{etale cover of mod stack}
    \coprod_{x \in \mM} [U_x / G_x] \lr \mM
\end{align}
where, by Remark~\ref{d-crit is local model}, each quotient stack $[U_x / G_x]$ is equipped with data of a local model 
\begin{align} \label{local model on cover}
    \Lambda_{[U_x/G_x]} = (U_x, V_x, G_x, F_{V_x}, \omega_{V_x}, D_{V_x}, \phi_{V_x}),
\end{align} with $F_{V_x} = \Omega_{V_x}$, $\omega_{V_x} = df_x$, $D_{V_x} = 0$ and $\phi_{V_x} = \sigma_{V_x}$.

By Proposition~\ref{comparison of d-crit quot}, for any two closed points $x, y \in \mM$ we have an \'{e}tale cover 
\begin{align} \label{etale cover of overlaps}
    \coprod_{z} [T_z / G_z] \lr [U_x / G_x] \times_{\mM} [U_y / G_y]
\end{align}
by d-critical quotient charts centered at closed points $z \in [U_x / G_x] \times_{\mM} [U_y / G_y]$, so that we have local model structures on $[T_z / G_z]$ given by the data \begin{align} \label{local model on overlap}
    \Lambda_{[T_z / G_z]} = (T_z, S_z, G_z, F_{S_z}, \omega_{S_z}, D_{S_z}, \phi_{S_z}),
\end{align}with $F_{S_z} = \Omega_{S_z}$, $\omega_{S_z} = df_z$, $D_{S_z} = 0$ and $\phi_{S_z} = \sigma_{S_z}$.

For any morphism $[T_z / G_z] \to [U_x / G_x]$ in \eqref{etale cover of overlaps}, Proposition~\ref{comparison of d-crit quot} also produces the following data:
\begin{enumerate}
\item A commutative $G_z$-equivariant diagram
\begin{align*}
    \xymatrix{
    T_z \ar[r] \ar[d] & S_z \ar[d]^-{\theta_{zx}} \\
    U_x \ar[r] & V_x
    }
\end{align*}
where the horizontal arrows are closed embeddings and $\theta_{zx}$ is unramified, inducing the \'{e}tale morphism $[T_z / G_z] \to [U_x / G_x]$ in \eqref{etale cover of overlaps}.
\item A surjective $G_z$-equivariant morphism $\eta_{zx} \colon F_{V_x} |_{S_z} \to F_{S_z}$, compatible with the morphisms $\phi_{V_x}|_{S_z}$ and $\phi_{S_z}$, such that $\omega_z' := \eta_{zx} (\omega_{V_x}|_{S_z})$ and $\omega_{S_z}$ are $\Omega$-equivalent (cf. \cite[Definition~5.9]{KLS}). 
\item An induced isomorphism 
$$\eta_{zx} \colon  h^1(\kK_{[U_x / G_x]})|_{[T_z / G_z]} \lr h^1(\kK_{[T_z / G_z]}). $$
\item Isomorphisms 
$$\psi_{xy} \colon h^1(\kK_{[U_x / G_x]})|_{[U_x / G_x] \times_{\mM} [U_y / G_y]} \lr h^1(\kK_{[U_y / G_y]})|_{[U_x / G_x] \times_{\mM} [U_y / G_y]}$$
satisfying the cocycle condition and $\psi_{xy}|_{[T_z / G_z]} = \eta_{zy}^{-1} \circ \eta_{zx}$.
\end{enumerate}

We abstract this situation by saying that the data of a strongly \'{e}tale cover~\eqref{etale cover of mod stack} equipped with a local model structure~\eqref{local model on cover}, an \'{e}tale cover~\eqref{etale cover of overlaps} with a local model structure~\eqref{local model on overlap} and compatibility data described by items (1)-(4) above constitute a \textbf{VC package} for an Artin stack $\mM$ satisfying the conditions of Theorem-Construction~\ref{Kirwan blowup of stack with GMS}. Note that if $\mM$ is an Artin stack of DT type then, by our discussion, we may take the local model structures to correspond to invariant d-critical quotient charts, however in general we only require local models that satisfy the compatibilities (1)-(4).
\end{setup} 

For a detailed account of $\Omega$-equivalence, we refer the reader to \cite[Section~5]{KLS}. For the purposes of our discussion, it suffices to note that $\omega_z'$ being $\Omega$-equivalent to $\omega_{S_z}$ essentially means that we have another local model structure
$$\Lambda_{[T_z / G_z]}' = (T_z, S_z, G_z, F_{S_z}, \omega_z', D_{S_z}, \phi_{S_z})$$
on $[T_z / G_z]$, such that the induced perfect obstruction theory $E_{[T_z^s / G_z]}'$ naturally has the same obstruction sheaf and gives the same obstruction assignment with $E_{[T_z^s / G_z]}$.

The data of Setup~\ref{VC package - setup of stack of DT type} are preserved under Kirwan blowups and hence lift to the intrinsic stabilizer reduction of a stack of DT type.

\begin{theorem} \label{intr stab red admits VC}
Let $\mM$ be an Artin stack of DT type. Then its intrinsic stabilizer reduction $\tilde{\mM}$ admits a VC package.
\end{theorem}

\begin{proof}
By the discussion in Setup~\ref{VC package - setup of stack of DT type}, $\mM$ admits a VC package where moreover all local models corresponds to d-critical quotient charts.

By the construction of $\tilde{\mM}$, it suffices to show that if a stack $\nN$ (satisfying the conditions of Theorem-Construction~\ref{Kirwan blowup of stack with GMS}) admits a VC package, then its Kirwan blowup $\hat{\nN}$ admits a canonically induced VC package.

By the construction in the proof of Proposition~\ref{comparison of d-crit quot}, Proposition~\ref{prop 4.9} and Theorem-Construction~\ref{Kirwan blowup of stack with GMS}, we may obtain new \'{e}tale covers~\ref{etale cover of mod stack} and \ref{etale cover of overlaps} for $\hat{\nN}$ by replacing every stack by its Kirwan blowup.

By Proposition~\ref{propagation of local model structure}, there are canonically induced local model structures on these \'{e}tale covers. It is moreover clear that condition (1) is satisfied.

Conditions (2)-(4) follow from \cite[Lemma~6.1]{KLS}. This concludes the proof.
\end{proof}

Now, a VC package on a Deligne--Mumford stack naturally provides data of a semi-perfect obstruction theory.

\begin{theorem} \label{VC on dm induces spot}
Let $\mM$ be a Deligne--Mumford stack equipped with a VC package. Then $\mM$ admits an induced semi-perfect obstruction theory and thus a virtual fundamental cycle $[\mM]\virt \in A_*(\mM)$.
\end{theorem}

\begin{proof}
Since $\mM$ is Deligne--Mumford, all quotient stacks in \eqref{etale cover of mod stack} are also Deligne--Mumford and satisfy $[U_x / G_x] = [U_x^s / G_x]$, since all closed points are now stable.

Thus, by Definition~\ref{pot on stable locus}, the local model structures $\Lambda_{[U_x / G_x]}$ induce perfect obstruction theories $$E_x := E_{[U_x / G_x]} \xrightarrow{\phi_x} L_{[U_x/G_x]}^{\geq -1}$$ on each $[U_x / G_x]$.

Conditions (3) and (4) of a VC package give descent data for the obstruction sheaves $h^1(E_x^\vee)$, while conditions (1) and (2) and the properties of $\Omega$-equivalence imply that these descent data are compatible with obstruction assignments of infinitesimal lifting problems.

Using Remark~\ref{spot schemes to dm stacks} and Theorem-Construction~\ref{spot virtual cycle def}, we obtain a semi-perfect obstruction theory on $\mM$ and a virtual fundamental cycle $[\mM]\virt \in A_*(\mM)$.
\end{proof}

Combining the above two theorems, we arrive at the main result of this section.

\begin{theorem} \label{dt type admits vir cyc}
Let $\mM$ be an Artin stack of DT type. Then the intrinsic stabilizer reduction $\tilde{\mM}$ admits a semi-perfect obstruction theory of virtual dimension zero and an induced, canonical virtual fundamental cycle $[\tilde{\mM}]\virt \in A_0(\tilde{\mM})$.
\end{theorem}

\begin{proof}
The existence of the semi-perfect obstruction theory follows immediately from Theorems~\ref{intr stab red admits VC} and \ref{VC on dm induces spot}.

The virtual dimension is zero, since for a d-critical chart the rank of $F_V = \Omega_V$ is equal to $\dim V$ and the virtual dimension $\dim V - \rk F_V$ is preserved under Kirwan blowups.

Finally, any two different \'{e}tale covers~\eqref{etale cover of mod stack} that are part of a VC package of $\mM$ can be refined to a common \'{e}tale cover. It is routine to check that the virtual fundamental cycle thus obtained is the same using this common refinement.
\end{proof}

\begin{remark}
The existence of the intrinsic stabilizer reduction of $\mM$ and its obstruction theory only uses its $(-1)$-shifted symplectic derived enhancement to deduce the existence of a d-critical structure on $\mM$. We could have thus used a weaker notion of stacks of DT type, replacing condition (3) in Definition~\ref{stack of dt type} by the requirement that $\mM$ admits a d-critical structure. However, we will be interested in replicating our arguments in the relative case where $\mM$ is a stack over a base smooth scheme $S$. In that context, there is no well-developed theory of d-critical stacks to the author's knowledge and the existence of a derived enhancement will allow us to still perform our constructions.
\end{remark}

\section{Donaldson--Thomas Invariants of Derived Objects} \label{sec DTK invariants}

In this section, $W$ denotes a smooth, projective Calabi--Yau threefold over $\bC$. We first describe the stability conditions $\sigma$ on $D^b(\Coh W)$ that we will be interested in. We then quote results in \cite{AlpHalpHein} which imply that moduli stacks of $\sigma$-semistable complexes are stacks of DT type and explain how to rigidify the $\bC^\ast$-scaling automorphisms of objects. Finally, we define generalized DT invariants via Kirwan blowups and show their deformation invariance.

\subsection{Stability conditions} By \cite{Lieblich}, there is an Artin stack $\pP := \Perf(W)$ of (universally gluable) perfect complexes on $W$, which is locally of finite type and has separated diagonal. Following \cite{AlpHalpHein}, we will consider the following type of stability condition.
\begin{definition} \emph{(Stability condition)} \label{stab cond} A stability condition $\sigma$ on $D^b(\Coh W)$ consists of the following data: 
\begin{enumerate}
    \item A heart $\aA \sub D^b ( \Coh W )$ of a t-structure on $D^b(\Coh W)$. Let $\pP_\aA$ denote the stack of perfect complexes in $\aA$.
    \item A vector $\gamma \in H^* (W, \bQ)$. Let $\pP_\aA^\gamma$ denote the stack of perfect complexes in $\aA$ with Chern character $\gamma$.   
    \item A locally constant function
\begin{align*}
    p_\gamma \colon \pi_0 (\pP_\aA) \lr V
\end{align*}
where $V$ is a totally ordered abelian group, $p_\gamma (E) = 0$ for $E \in \pP_\aA^\gamma$ and $p_\gamma$ is additive so that $p_\gamma(E \oplus F) = p_\gamma(E) + p_\gamma(F)$.
\end{enumerate}

We say that $E \in \pP_\aA^\gamma$ is semistable if for any subobject $F \subseteq E$ we have $p_\gamma(F) \leq 0$ and stable if $p_\gamma (F) < 0$. If $E$ is not semistable, we say it is unstable.
\end{definition}

In order for the stack of semistable objects to be of DT type, we will need to consider stability conditions satisfying certain properties.

\begin{definition} \emph{(Nice stability condition)} \label{nice stab cond} Given a stability condition $\sigma$ on $D^b(\Coh W)$, let $\mM$ be the stack of $\sigma$-semistable objects in $\pP_\aA^\gamma$. We say that $\sigma$ is nice if the following hold:
\begin{enumerate}
\item $\mM$ is an Artin stack of finite type.
\item $\mM$ is an open substack of $\Perf(W)$.
\item $\mM$ satisfies the existence part of the valuative criterion of properness. We then say that $\mM$ is quasi-proper or universally closed.
\end{enumerate}
\end{definition}

\begin{remark} The following are examples of nice stability conditions, by the results of the mentioned authors:
\begin{enumerate}
\item A Bridgeland stability condition in the sense of Piyaratne--Toda \cite{TodaPiya} and Li \cite{quinticstab}.
\item A polynomial stability condition in the sense of Lo \cite{Lo,Lo2}.
\item Gieseker and slope stability. These are examples of a weak stability condition in the sense of Joyce--Song \cite[Definition~3.5]{JoyceSong}, where we take $\aA = \Coh W$ and $K(A) = N(W)$.
\end{enumerate}
\end{remark}

\subsection{Moduli stacks of semistable complexes are of DT type} The following theorem is an application of \cite[Theorem~7.25]{AlpHalpHein} to our context.

\begin{theorem} \label{6.13}
Let $\sigma$ be a nice stability condition on $D^b(\Coh W)$ and let $\mM$ denote the stack of $\sigma$-semistable complexes. Then $\mM$ is an Artin stack of DT type with proper good moduli space $M$ and affine diagonal.
\end{theorem}

\begin{proof}
By the niceness of $\sigma$, $\mM$ is of finite type. By the results of \cite{PTVV}, $\Perf(W)$ is the truncation of a $(-1)$-shifted symplectic derived Artin stack. Since $\mM$ is an open substack of $\Perf(W)$, it is also such a truncation. By \cite[Theorem~7.25]{AlpHalpHein}, $\mM$ admits a good moduli space $\pi \colon \mM \to M$ such that $M$ is separated. By Proposition~\ref{properties of good mod spaces}, $\pi$ is universally closed and since $\mM$ is universally closed as well, it follows that $M$ is universally closed. Thus $M$ must be proper. 

$\mM$ has affine diagonal by \cite[Lemma 7.19]{AlpHalpHein}.
\end{proof}

\subsection{$\bC^\ast$-rigidified intrinsic stabilizer reduction} \label{Subsection 5.3} Let $\mM$ be as in Theorem~\ref{6.13}. We denote $T = \bC^\ast$. 

To obtain a meaningful nonzero DT invariant, it will be necessary to rigidify the $\bC^\ast$-scaling automorphisms of objects in $\mM$. 

For each family of complexes $E_S \in \mM(S)$ there exists an embedding
$$ \bG_m (S) \to \Aut(E_S)$$
which is compatible with pullbacks and moreover $\bG_m (S)$ is central. In the terminology used in \cite{AGV}, we say that $\mM$ has a $\bG_m$-2-structure.

Using the results of \cite{AOV} or \cite{AGV}, we may take the $\bG_m$-rigification $\mM \!\!\fatslash \bG_m$ of $\mM$. From the properties of rigification, for any point $x \in {\mM}$, one has $\Aut_{{\mM}\!\!\fatslash \bG_m}(x) = \Aut_{ {\mM}}(x) / T$. In particular, if $x \in \mM^s$ is stable, then $\Aut_{\mM\!\!\fatslash \bG_m}(x) = \lbrace \id \rbrace$. 

Even though $\mM\!\!\fatslash \bG_m$ is not a stack of DT type, we show in the next two propositions that we can construct its intrinsic stabilizer reduction and equip it with a semi-perfect obstruction theory of virtual dimension zero.

\begin{proposition} \label{prop 6.5}
$\mM\!\!\fatslash \bG_m$ satisfies the conditions of Theorem~\ref{intrinsic stab red of stack with gms}. It thus admits an intrinsic stabilizer reduction, which we denote by $\tilde{\mM}^{\bC^*}$.
\end{proposition}

\begin{proof}
By \cite[Theorem~5.1]{Romagny} or \cite[Theorem~C.1.1]{AOV}, $\mM\!\!\fatslash \bG_m$ has the same good moduli space $M$. Since the morphism $\mM \to \mM\!\!\fatslash \bG_m$ is a $\bG_m$-gerbe, it follows that $\mM\!\!\fatslash \bG_m$ and the good moduli space morphism $\mM\!\!\fatslash \bG_m \to M$ have affine diagonal.
\end{proof}

\begin{proposition}
$\mM\!\!\fatslash \bG_m$ is equipped with a VC package.
\end{proposition}

\begin{proof}
$\mM$ is a stack of DT type and hence admits a VC package by Setup~\ref{VC package - setup of stack of DT type}. By the existence of the $\bG_m$-2-structure, $T$ naturally embeds in each stabilizer group $G_x$. Replacing all such groups by their quotients $G_x / T$, we obtain a VC package for $\mM\!\!\fatslash \bG_m$.
\end{proof}

\begin{remark}
In the case of semistable sheaves treated in \cite{KLS}, rigidification is much simpler, since the moduli stack is a global GIT quotient stack $\mM = [ X / G ]$ where $G = \GL (N, \bC)$, and then one may work with $[ X / \mathrm{PGL} (N, \bC) ]$ as the $\bG_m$-rigidication.
\end{remark}

We thus obtain the following.

\begin{thm-defi} \label{rigidified intr stab red}
$\tilde{\mM}^{\bC^*}$ is called the $\bC^\ast$-rigidified intrinsic stabilizer reduction of $\mM$. It is a proper DM stack with a semi-perfect obstruction theory of virtual dimension zero.
\end{thm-defi}

\begin{proof}
For properness, the good moduli space $q \colon \tilde{\mM}^{\bC^\ast} \to \tilde{M}^{\bC^\ast}$ is proper. Since $M$ is proper and $\tilde{M}^{\bC^\ast}$ is proper over $M$, $\tilde{\mM}^{\bC^\ast}$ is proper.

Everything else follows from Theorems~\ref{intr stab red admits VC}, \ref{VC on dm induces spot} as applied in the proof of Theorem~\ref{dt type admits vir cyc}.
\end{proof}

\begin{remark}
By identical reasoning, all of the above hold in greater generality when $\mM$ is an Artin stack of DT type with a $\bG_m$-2-structure.
\end{remark}

\subsection{Generalized DT invariants via Kirwan blowups} Suppose as before that $\mM$ is an Artin stack of DT type parametrizing $\sigma$-semistable objects in a heart $\aA$ of a $t$-structure on $D^b(\Coh W)$ with Chern character $\gamma$, for a nice stability condition $\sigma$. Then, by Theorem-Definition~\ref{rigidified intr stab red}, there is an induced $\bC^\ast$-rigidified intrinsic stabilizer reduction $\tilde{\mM}^{\bC^\ast}$ with a good moduli space $\tilde{M}^{\bC^\ast}$ and an induced semi-perfect obstruction theory and virtual cycle of dimension zero. We can now state the main theorem of this paper.

\begin{thm-defi} \label{DTK invariant of complexes}
Let $W$ be a smooth, projective Calabi--Yau threefold, $\sigma$ a nice stability condition on a heart $\aA \subset D^b (\Coh W)$ of a $t$-structure, as in Definition~\ref{stab cond}, $\gamma \in H^*(W)$ and let $\mM$ denote the stack parametrizing $\sigma$-semistable complexes with Chern character $\gamma$. Then we may define the associated generalized Donaldson--Thomas invariant via Kirwan blowups (DTK invariant) as $$\mathrm{DTK}\left( \mM \right) := \deg \ [ \tilde{\mM}^{\bC^\ast} ]^{\mathrm{vir}} \in \bQ.$$
\end{thm-defi}

\begin{remark}
In \cite{KiemSavvas}, the results of \cite{KLS} and the present paper are refined to define a virtual structure sheaf $[\oO_{\tilde{\mM}}\virt] \in K_0(\tilde{\mM}^{\bC^\ast})$ and a corresponding $K$-theoretic generalized DTK invariant.
\end{remark}

\subsection{Relative theory and deformation invariance} Let $C$ be a smooth quasi-projective scheme over $\bC$ and $W \to C$ a smooth, projective family of Calabi--Yau threefolds. Without loss of generality, we assume that $H^*(W_t, \bQ)$ stays constant for $t \in C$ and identify it with $H^*(W_0, \bQ)$ where $W_0$ is the fiber of the family over a point $0 \in C$. Let $\gamma \in H^*(W_0, \bQ)$. Moreover, let $\Perf(W/S)$ denote the stack of (universally gluable) perfect complexes on the morphism $W \to C$ as in \cite{Lieblich}.

We consider families $\sigma$ of stability conditions $\sigma_t$, where, for each $t \in C$, $\sigma_t$ is a stability condition on $D^b(\Coh W_t)$ as in Definition~\ref{nice stab cond} with Chern character $\gamma \in H^*(W_t, \bQ) = H^*(W_0, \bQ)$. Let $\mM \to C$ be the stack parametrizing relatively semistable objects in $D^b(\Coh W)$, i.e. perfect complexes $E$ such that the derived restriction $E_t := E|_{W_t}$ is $\sigma_t$-semistable for all $t \in C$. We require that the conditions characterizing a nice stability condition hold relative to the base $C$ as follows.

\begin{setup} \label{relative setup for M over S}
We say that the family of stability conditions $\sigma_t$ is nice if the following hold.
\begin{enumerate}
    \item $\mM$ is an Artin stack of finite type.
    \item $\mM$ is an open substack of $\Perf(W/C)$.
    \item $\mM \to C$ admits a good moduli space $M \to C$, proper over $C$.
\end{enumerate}
\end{setup}

\begin{remark}
When we have a GIT description of $\mM \to C$, then the above conditions are satisfied. This is the case for Gieseker and slope stability of coherent sheaves.
\end{remark}

When we are in the situation of the above setup, all of our results extend to the relative case. 

\begin{theorem} \label{def inv of dtk}
Let $W \to C$ be a smooth, projective family of Calabi--Yau threefolds over a smooth, quasi-projective scheme $C$, $\lbrace \sigma_t \rbrace_{t \in C}$ a nice family of stability conditions on $D^b( \Coh W)$, $\gamma \in H^*(W_0, \bQ)$ and $\mM \to C$ the stack of fiberwise $\sigma_t$-semistable objects of Chern character $\gamma$ in $D^b(\Coh W)$.

Then there exists an induced $\bC^\ast$-rigidified intrinsic stabilizer reduction $\tilde{\mM} \to C$, a DM stack proper over $C$, endowed with a semi-perfect obstruction theory  and a virtual fundamental cycle $[\tilde{\mM}]^{\mathrm{vir}} \in A_0(\tilde{\mM})$. 

Moreover, the fiber $\tilde{\mM}_t$ over $t \in C$ is the $\bC^\ast$-rigidified intrinsic stabilizer reduction of $\mM_t$ and the obstruction theory pulls back to the one constructed in the absolute case, so that if $i_t \colon \tilde{\mM}_t \to \tilde{\mM}$ is the inclusion, we have
$$i_t^![\tilde{\mM}]\virt = [\tilde{\mM}_t]\virt \in A_0(\tilde{\mM}_t).$$
\end{theorem}

\begin{proof}
This is the generalization of \cite[Theorem~7.17]{KLS} adapted to our context. We briefly explain the steps and necessary changes.

$\mM$ admits a $\bG_m$-2-structure. By conditions (1) and (3) of Setup~\ref{relative setup for M over S}, $\mM\!\!\fatslash \bG_m$ satisfies the conditions of Theorem~\ref{intrinsic stab red of stack with gms}, so using Proposition~\ref{prop 6.5} we obtain the  $\bC^\ast$-rigidified intrinsic stabilizer reduction $\tilde{\mM}^{\bC^\ast} \to \mM \to C$.

To see that the fiber over $t \in C$ is $\tilde{\mM}_t^{\bC^\ast}$, observe that by definition $\tilde{\mM}^{\bC^\ast}$ is obtained by taking iterated Kirwan blowups of a cover by strongly \'{e}tale quotient charts (where $T = \bC^\ast$)
\begin{align*}
    \xymatrix{
    \nN_x := [U_x / (G_x / T)] \ar[d] \ar[r] & \mM\!\!\fatslash \bG_m \ar[d] \\
    N_x := U_x \git G_x \ar[r] & M \ar[d] \\
     & C.
    }
\end{align*}

Taking fibers over $t \in C$, we obtain a cover of $\mM_t \!\!\fatslash \bG_m$ by strongly \'{e}tale quotient charts, so it suffices to show that $(\hat{\nN}_x)_t = \hat{(\nN_x)}_t$. But this is true by the same argument used in the proof of \cite[Proposition~7.4]{KLS}.

Constructing the obstruction theory of $\tilde{\mM}^{\bC^\ast}$ is slightly more subtle. First, we observe that there is an obvious generalization of Setup-Definition~\ref{local model} to define local model structures on $G$-schemes $U$ over the base scheme $C$, cf. \cite[Definition~7.11]{KLS}. 

Using the fact that the morphism $\mM \to C$ is the truncation of a $(-1)$-shifted symplectic derived stack over $C$, we can show that $\mM\!\!\fatslash \bG_m$ admits a VC package, as in Setup~\ref{VC package - setup of stack of DT type}. The only difference is that now in \eqref{local model on cover}, we have $F_{V_x} = \Omega_{V_x/C}$ and $\omega_{V_x}$ is a $G_x$-invariant $1$-form, and similarly in \eqref{local model on overlap}, $F_{S_z} = \Omega_{S_z / C}$ and $\omega_{S_z}$ is a $G_z$-invariant $1$-form. The construction of the VC package follows verbatim the reasoning of Subsections~7.2 and 7.3 in \cite{KLS}. We refer the reader there for details.

Finally, to see that the restriction of the semi-perfect obstruction theory of $\tilde{\mM}^{\bC^\ast}$ to $\tilde{\mM}^{\bC^\ast}_t$ agrees with the absolute semi-perfect obstruction theory constructed using the d-critical structure of $\mM_t \!\!\fatslash \bG_m$, we use the fact that for any choice of relative local models in the VC package of $\mM\!\!\fatslash \bG_m$, the sections $\omega_{V_x}$ can be taken so that $\omega_{V_x}|_{(V_x)_t}$ is $\Omega$-equivalent to an exact $1$-form $df_x$ induced by the d-critical structure of $\mM_t \!\!\fatslash \bG_m$. As a consequence of properties of $\Omega$-equivalence (cf. \cite[Lemma~5.14]{KLS}), the construction of the semi-perfect obstruction theory of $\tilde{\mM}^{\bC^\ast}_t$ can be performed equivalently using $\omega_{V_x}|_{(V_x)_t}$ instead of $df_x$ in the VC package of $\mM_t \!\!\fatslash \bG_m$. The necessary arguments are carried out in detail in \cite[Subsection~7.4]{KLS}.
\end{proof}

In the case of Bridgeland stability conditions constructed in \cite{TodaPiya} and \cite{quinticstab}, the results of \cite{familystab} imply that we get nice families of Bridgeland stability conditions.

As an immediate corollary, we have the following theorem.

\begin{theorem}
The generalized DT invariant via Kirwan blowups for $\sigma$-semistable objects on Calabi--Yau threefolds, where $\sigma$ is a Bridgeland stability condition as in \cite{TodaPiya} and \cite{quinticstab}, is invariant under deformations of the complex structure of the Calabi--Yau threefold.
\end{theorem}

\begin{remark}
In the case of Gieseker stability and slope stability of coherent sheaves, deformation invariance follows directly from the results of \cite{KLS}.
\end{remark}

\bibliography{Master}
\bibliographystyle{alpha}

\end{document}